\documentclass[12pt, a4paper]{article}

\usepackage{amssymb}
\usepackage{amsmath}

  \newtheorem{theorem}{Theorem}
 \newtheorem{definition}{Definition}
  \newtheorem{corollary}{Corollary}
    \newtheorem{lemma}{Lemma}

    \newcommand\vu{{\bf u}}
    \newcommand\vv{{\bf  v}}
    \newcommand\vw{{\bf  w}}
    \newcommand\vb{{\bf b}}
    \newcommand\vN{{\bf  \nabla}}
    
    \newcommand\Proof{\noindent {\bf Proof : }}
    \newcommand \Endproof{\hfill $\diamond$}

\begin{document}

\title{Weak solutions for Navier--Stokes equations with
initial data in   weighted $L^2$ spaces.}
\author{Pedro Gabriel Fern\'andez-Dalgo\footnote{LaMME, Univ Evry, CNRS, Universit\'e Paris-Saclay, 91025, Evry, France } \footnote{e-mail : pedro.fernandez@univ-evry.fr, pedrodalgo16@gmail.com} ,  Pierre Gilles Lemari\'e--Rieusset\footnote{LaMME, Univ Evry, CNRS, Universit\'e Paris-Saclay, 91025, Evry, France} \footnote{e-mail : pierregilles.lemarierieusset@univ-evry.fr}}
\date{}\maketitle

\begin{abstract}
We show the existence of global weak solutions of the 3D Navier-Stokes equations with initial velocity
in the weighted spaces  $L^2_{w_\gamma}$, where $w_\gamma(x)=(1+\vert x\vert)^{-\gamma}$ and $0<\gamma\leq 2$, using new energy controls. As
application we give a new proof of the existence of global weak discretely self-similar solutions of the 3D
Navier--Stokes equations for discretely self-similar initial velocities which are locally square integrable.\end{abstract}
 


\section{Introduction.}

Infinite-energy weak Leray solutions  to the Navier--Stokes equations were introduced by  Lemari\'e-Rieusset
in 1999  \cite{LR99} (they are presented more completely in \cite{LR02}  and  \cite{LR16}). This has allowed to show the
existence of local weak solutions for a  uniformly locally square integrable initial data.

Other constructions of  infinite-energy solutions for locally uniformly square integrable
initial data were given in 2006 by Basson \cite{Ba06} and in 2007 by Kikuchi and Seregin \cite{KS07}. These
solutions allowed   Jia and Sverak \cite{JS14} to construct in 2014 the self-similar solutions for large
(homogeneous of degree -1) smooth data. Their result has been extended in 2016 by  Lemari\'e-Rieusset
\cite{LR16} to solutions for rough locally square integrable data. We remark that an homogeneous (of degree
-1) and locally square integrable data is automatically uniformly locally $L^2$.

Recently,  Bradshaw and Tsai \cite{BT19}  and Chae and Wolf \cite{CW18}  considered the case of  solutions which are self-similar  according  to a discrete subgroup
of dilations. Those solutions are related to an initial data which is self-similar  only for a discrete group  of dilations; in contrast to the case of self-similar solutions for all dilations,  such an initial data, when  
 locally  $L^2$,  is not necessarily uniformly locally  $L^2$, therefore their
results are no  consequence of constructions described by Lemari\'e-Rieusset in \cite{LR16}.

  In this paper, we  construct
an alternative theory to obtain   infinite-energy global weak solutions for large initial data, which
include the discretely self-similar locally square integrable data. More specifically, we consider the weights
  $$ w_\gamma(x)=\frac 1{(1+\vert x\vert)^\gamma}$$ with $0<\gamma$, and the spaces
  $$ L^2_{w_\gamma}=L^2(w_\gamma \, dx).$$
  Our main theorem is the following one :

 \begin{theorem} \label{weightedNS} Let $0<\gamma\leq 2$. If $\vu_{0}$ is a divergence-free vector field such that $\vu_0\in L^2_{w_\gamma}(\mathbb{R}^3)$ and if $\mathbb{F}$ is a tensor $\mathbb{F}(t,x)=\left(F_{i,j}(t,x)\right)_{1\leq i,j\leq 3}$ such that $\mathbb{F}\in L^2((0,+\infty), L^2_{w_\gamma})$, then the Navier--Stokes equations with initial value $\vu_0$
 \begin{equation*}  (NS) \left\{ \begin{matrix} \partial_t \vu= \Delta \vu  -(\vu\cdot \vN)\vu- \vN p +\vN\cdot \mathbb{F} \cr \cr \vN\cdot \vu=0,    \phantom{space space} \vu(0,.)=\vu_0
 \end{matrix}\right.\end{equation*}
 has a global weak solution $\vu$ such that :
 \begin{itemize} 
 \item[$\bullet$] for every $0<T<+\infty$, $\vu$ belongs to $L^\infty((0,T), L^2_{w_\gamma})$ and $\vN\vu$ belongs to $L^2((0,T),L^2_{w_\gamma})$
 \item[$\bullet$] the pressure $p$  is related to $\vu$ and $\mathbb{F}$ through the Riesz transforms $R_i =\frac{\partial_i}{\sqrt{-\Delta}}$ by the formula
 $$ p =\sum_{i=1}^3\sum_{j=1}^3 R_iR_j(u_iu_j-F_{i,j})$$
 where, for every $0<T<+\infty$,   $\sum_{i=1}^3\sum_{j=1}^3 R_iR_j(u_iu_j)$  belongs to $L^{4}((0,T),L^{6/5}_{w_{\frac {6\gamma}5}})$ and  $\sum_{i=1}^3\sum_{j=1}^3 R_iR_j F_{i,j}  $ belongs to $L^{2}((0,T),L^{2}_{w_\gamma})$
 \item[$\bullet$] the map $t\in [0,+\infty)\mapsto \vu(t,.)$ is weakly continuous from $[0,+\infty)$ to $L^2_{w_\gamma}$, and is strongly continuous at $t=0$ :
 $$ \lim_{t\rightarrow 0} \|\vu(t,.)-\vu_0\|_{L^2_{w_\gamma}}=0.$$
 \item[$\bullet$] the solution $\vu$ is suitable : there exists a non-negative locally finite measure $\mu$ on $(0,+\infty)\times\mathbb{R}^3$ such that
 $$ \partial_t(\frac {\vert\vu\vert^2}2)=\Delta(\frac {\vert\vu\vert^2}2)-\vert\vN \vu\vert^2- \vN\cdot\left( (\frac{\vert\vu\vert^2}2+p)\vu\right) + \vu\cdot(\vN\cdot\mathbb{F})-\mu.$$
 \end{itemize}
 
 In particular, we have the energy controls
 \begin{equation*}\begin{split}
 \|\vu(t,. )\|_{L^2_{w_\gamma}}^2 &+2\int_0^t \|\vN \vu(s,.)\|_{L^2_{w_\gamma}}^2\, ds
 \\  \leq& \|\vu_0 \|_{L^2_{w_\gamma}}^2 -\int_0^t \int \vN\vert\vu\vert^2\cdot\vN w_\gamma \, dx \, ds  +\int_0^t\int (\vert\vu\vert^2+2p) \vu\cdot \vN(w_\gamma)\,  \, dx \, ds \\ &\phantom{space space} -2\sum_{i=1}^3\sum_{j=1}^3 \int_0^t\int F_{i,j} (\partial_i u_j) w_\gamma +F_{i,j} u_i \partial_j(w_\gamma)  \, dx \, ds
 \end{split}\end{equation*}
 and
 \begin{equation*}   \|\vu(t,. )\|_{L^2_{w_\gamma}}^2  \leq \|\vu_0 \|_{L^2_{w_\gamma}}^2  +C_\gamma\int_0^t \|  \mathbb{F}(s,.)\|^2_{L^2_{w_\gamma}}\, ds + C_\gamma  \int_0^t   \|\vu(s,. )\|_{L^2_{w_\gamma}}^2  + \|\vu(s,. )\|_{L^2_{w_\gamma}}^6  \, ds\end{equation*}
 \end{theorem}
 
 A  key tool for proving Theorem \ref{weightedNS} and for applying it to the study of   discretely self-similar solutions is given by  the following a priori estimates for an advection-diffusion problem :
 \begin{theorem}\label{estimates}
  Let $0<\gamma\leq 2$. Let $0<T<+\infty$. Let  $\vu_{0}$ be a divergence-free vector field such that $\vu_0\in L^2_{w_\gamma}(\mathbb{R}^3)$ and   $\mathbb{F}$ be  a tensor $\mathbb{F}(t,x)=\left(F_{i,j}(t,x)\right)_{1\leq i,j\leq 3}$ such that $\mathbb{F}\in L^2((0,T), L^2_{w_\gamma})$.  Let $\vb$ be a time-dependent divergence free vector-field ($\vN\cdot\vb=0$) such that $\vb\in L^3((0,T),L^3_{w_{3\gamma/2}})$.
  
 Let $\vu$ be a solution of the following advection-diffusion problem 
 \begin{equation*}  (AD) \left\{ \begin{matrix} \partial_t \vu= \Delta \vu  -(\vb\cdot \vN)\vu- \vN p +\vN\cdot \mathbb{F} \cr \cr \vN\cdot \vu=0,    \phantom{space space} \vu(0,.)=\vu_0
 \end{matrix}\right.\end{equation*}
be such that :
 \begin{itemize} 
 \item[$\bullet$]  $\vu$ belongs to $L^\infty((0,T), L^2_{w_\gamma})$ and $\vN\vu$ belongs to $L^2((0,T),L^2_{w_\gamma})$
 \item[$\bullet$] the pressure $p$  is related to $\vu$, $\vb$  and $\mathbb{F}$ through the Riesz transforms $R_i =\frac{\partial_i}{\sqrt{-\Delta}}$ by the formula
 $$ p=\sum_{i=1}^3\sum_{j=1}^3 R_iR_j(b_iu_j-F_{i,j})$$
 where   $\sum_{i=1}^3\sum_{j=1}^3 R_iR_j(b_iu_j)$  belongs to $L^{3}((0,T),L^{6/5}_{w_{\frac {6\gamma}5}})$ and  $\sum_{i=1}^3\sum_{j=1}^3 R_iR_j F_{i,j}  $ belongs to $L^{2}((0,T),L^{2}_{w_\gamma})$
 \item[$\bullet$] the map $t\in [0,T)\mapsto \vu(t,.)$ is weakly continuous from $[0,T)$ to $L^2_{w_\gamma}$, and is strongly continuous at $t=0$ :
 $$ \lim_{t\rightarrow 0} \|\vu(t,.)-\vu_0\|_{L^2_{w_\gamma}}=0.$$
 \item[$\bullet$]  there exists a non-negative locally finite measure  $\mu$ on $(0,T)\times\mathbb{R}^3$ such that
 \begin{equation}\label{energloc} \partial_t(\frac {\vert\vu\vert^2}2)=\Delta(\frac {\vert\vu\vert^2}2)-\vert\vN \vu\vert^2- \vN\cdot\left( \frac{\vert\vu\vert^2}2\vb\right)-\vN\cdot(p\vu) + \vu\cdot(\vN\cdot\mathbb{F})-\mu.\end{equation}
 \end{itemize}
 
Then, we have the energy controls
 \begin{equation*}\begin{split}
 \|\vu(t,. )\|_{L^2_{w_\gamma}}^2 &+2\int_0^t \|\vN \vu(s,.)\|_{L^2_{w_\gamma}}^2\, ds
 \\  \leq& \|\vu_0 \|_{L^2_{w_\gamma}}^2  -\int_0^t \int   \vN\vert\vu\vert^2\cdot\vN w_\gamma\,   dx \, ds +\int_0^t\int \vert\vu\vert^2 \vb\cdot \vN(w_\gamma) \, dx \, ds \\ &+ 2\int_0^t\int  p\vu\cdot \vN(w_\gamma)\, dx\, ds  -2\sum_{i=1}^3\sum_{j=1}^3 \int_0^t\int F_{i,j} (\partial_i u_j) w_\gamma +F_{i,j} u_i \partial_j(w_\gamma)\, dx\, ds
 \end{split}\end{equation*}
 and 
 \begin{equation*} \begin{split}  \|\vu(t,. )\|_{L^2_{w_\gamma}}^2 &+  \int_0^t  \|\vN\vu\|_{L^2_{w_\gamma}}^2\, ds  \\ \leq& \|\vu_0 \|_{L^2_{w_\gamma}}^2  +C_\gamma \int_0^t \|  \mathbb{F}(s,.)\|^2_{L^2_{w_\gamma}}\, ds + C_\gamma  \int_0^t  (1+  \|\vb(s,. )\|_{L^3_{w_{3\gamma/2}}}^2 )  \|\vu(s,. )\|_{L^2_{w_\gamma}}^2  \, ds\end{split}\end{equation*} where $C_\gamma$ depends only on $\gamma$ (and not on $T$, and not on $\vb$, $\vu$, $\vu_0$ nor $\mathbb{F}$).
 \end{theorem}
 
 In particular, we shall prove the following stability result :

 \begin{theorem}\label{stability}
  Let $0<\gamma\leq 2$. Let $0<T<+\infty$. Let  $\vu_{0,n}$ be  divergence-free vector fields such that $\vu_{0,n}\in L^2_{w_\gamma}(\mathbb{R}^3)$ and  $\mathbb{F}_n$ be    tensors  such that $\mathbb{F}_n\in L^2((0,T), L^2_{w_\gamma})$.  Let $\vb_n$ be   time-dependent divergence free vector-fields such that $\vb_n\in L^3((0,T),L^3_{w_{3\gamma/2}})$.
  
 Let $\vu_n$ be solutions of the following advection-diffusion problems 
 \begin{equation*}  (AD_n) \left\{ \begin{matrix} \partial_t \vu_n= \Delta \vu_n  -(\vb_n\cdot \vN)\vu_n- \vN p_n +\vN\cdot \mathbb{F}_n \cr \cr \vN\cdot \vu_n=0,    \phantom{space space} \vu_n(0,.)=\vu_{0,n}
 \end{matrix}\right.\end{equation*}
  such that :
 \begin{itemize} 
 \item[$\bullet$] $\vu_n$ belongs to $L^\infty((0,T), L^2_{w_\gamma})$ and $\vN\vu_n$ belongs to $L^2((0,T),L^2_{w_\gamma})$
 \item[$\bullet$] the pressure $p_n$  is related to $\vu_n$, $\vb_n$  and $\mathbb{F}_n$   by the formula
 $$ p_n=\sum_{i=1}^3\sum_{j=1}^3 R_iR_j(b_{n,i}u_{n,j}-F_{n,i,j})$$
 \item[$\bullet$] the map $t\in [0,T)\mapsto \vu_n(t,.)$ is weakly continuous from $[0,T)$ to $L^2_{w_\gamma}$, and is strongly continuous at $t=0$ :
 $$ \lim_{t\rightarrow 0} \|\vu_n(t,.)-\vu_{0,n}\|_{L^2_{w_\gamma}}=0.$$
 \item[$\bullet$]  there exists a non-negative locally finite measure $\mu_n$ on $(0,T)\times\mathbb{R}^3$ such that
 $$ \partial_t(\frac {\vert\vu_n\vert^2}2)=\Delta(\frac {\vert\vu_n\vert^2}2)-\vert\vN \vu_n\vert^2- \vN\cdot\left( \frac{\vert\vu_n\vert^2}2\vb_n\right)-\vN\cdot(p_n\vu_n) + \vu_n\cdot(\vN\cdot\mathbb{F}_n)-\mu_n.$$
 \end{itemize}
 
 If $\vu_{0,n}$ is strongly convergent to $\vu_{0,\infty}$ in $L^2_{w_\gamma}$,  if the sequence $\mathbb{F}_n$ is strongly convergent to $\mathbb{F}_\infty$  in $L^2((0,T), L^2_{w_\gamma})$,  and   if the sequence $\vb_n$ is bounded in $L^3((0,T),   L^3_{w_{3\gamma/2}})$, then there exists $p_\infty$, $\vu_\infty$,    $\vb_\infty$   and an increasing  sequence $(n_k)_{k\in\mathbb{N}}$ with values in $\mathbb{N}$ such that
  \begin{itemize} 
 \item[$\bullet$] $\vu_{n_k}$ converges *-weakly to $\vu_\infty$ in $L^\infty((0,T), L^2_{w_\gamma})$, $\vN\vu_{n_k}$ converges weakly to $\vN\vu_\infty$ in $L^2((0,T),L^2_{w_\gamma})$
 \item[$\bullet$] $\vb_{n_k}$ converges weakly to $\vb_\infty$ in $L^3((0,T), L^3_{w_{3\gamma/2}})$, $p_{n_k}$ converges weakly to $p_\infty$ in $L^{3}((0,T),L^{6/5}_{w_{\frac {6\gamma}5}})+L^{2}((0,T),L^{2}_{w_\gamma})$
  \item[$\bullet$] $\vu_{n_k}$ converges strongly  to $\vu_\infty$ in  $L^2_{\rm loc}([0,T)\times\mathbb{R}^3)$ : for every  $T_0\in (0,T)$ and every $R>0$, we have
  $$\lim_{k\rightarrow +\infty} \int_0^{T_0} \int_{\vert y\vert<R} \vert \vu_{n_k}(s,y)-\vu_\infty(s,y)\vert^2\, ds\, dy=0.$$
  \end{itemize}

 Moreover, $\vu_\infty$ is  a solution of the   advection-diffusion problem 
 \begin{equation*}  (AD_\infty) \left\{ \begin{matrix} \partial_t \vu_\infty= \Delta \vu_\infty  -(\vb_\infty\cdot \vN)\vu_\infty- \vN p_\infty +\vN\cdot \mathbb{F}_\infty \cr \cr \vN\cdot \vu_\infty=0,    \phantom{space space} \vu_\infty(0,.)=\vu_{0,\infty}
 \end{matrix}\right.\end{equation*}
 and is
  such that :
 \begin{itemize}  
 \item[$\bullet$] the map $t\in [0,T)\mapsto \vu_\infty(t,.)$ is weakly continuous from $[0,T)$ to $L^2_{w_\gamma}$, and is strongly continuous at $t=0$ :
 $$ \lim_{t\rightarrow 0} \|\vu_\infty(t,.)-\vu_{0,\infty}\|_{L^2_{w_\gamma}}=0.$$
 \item[$\bullet$]  there exists a non-negative locally finite measure $\mu_\infty$ on $(0,T)\times\mathbb{R}^3$ such that
 $$ \partial_t(\frac {\vert\vu_\infty\vert^2}2)=\Delta(\frac {\vert\vu_\infty\vert^2}2)-\vert\vN \vu_\infty\vert^2- \vN\cdot\left( \frac{\vert\vu_\infty\vert^2}2\vb_\infty\right)-\vN\cdot(p_\infty\vu_\infty) + \vu_\infty\cdot(\vN\cdot\mathbb{F}_\infty)-\mu_\infty.$$
 \end{itemize}
 \end{theorem}
 \section*{Notations.}
 All along the text, $C_\gamma$  is a positive constant whose value may change from line to line but which depends only on $\gamma$.
\section{The weights $w_\delta$.}
 We consider the weights  $w_\delta=\frac 1{(1+\vert x\vert)^\delta}$ where $0<\delta$ and $x\in\mathbb{R}^3$. A very important feature of those weights is the control of their gradients :
 \begin{equation} \vert\vN w_\delta(x)\vert =\delta \frac{w_\delta(x)}{1+\vert x\vert}\end{equation}

 \begin{lemma}[Muckenhoupt weights]\label{Lemmuck} If $0<\delta<3$ and $1<p<+\infty$, then $w_\delta$ belongs to the Muckenhoupt class $\mathcal{A}_p$.
 \end{lemma}
 
 \Proof We recall that a weight $w$ belongs to $\mathcal{A}_p(\mathbb{R}^3)$ for $1<p<+\infty$ if and only if it satisfies the reverse H\"older  inequality
\begin{equation}\label{muck} \!\!\! \sup_{x\in\mathbb{R}^3, R>0} \!\! \left(\frac 1{\vert B(x,R)\vert}  \!\int_{B(x,R)}\! \!\!\!\!w(y)\, dy\right)^{\frac 1 p}\!\! \left( \frac 1{\vert B(x,R)\vert}\! \int_{B(x,R)} \!\frac{dy}{w(y)^{\frac 1{p-1}}}\right)^{1-\frac 1 p} \!\! \!\!<+\infty.\end{equation}
 For all  $0<R\leq 1$ the inequality $\vert x-y\vert<R$ implies 
 $ \frac 1 2 (1+\vert x\vert) \leq 1+\vert y\vert\leq 2 (1+\vert x\vert)$, thus we can control the left
side in (\ref{muck}) for $w_\delta$  by $4^{\frac \delta p}$.
 
 For all $R > 1$  and $\vert x\vert>10 R$, we have that  the inequality $\vert x-y\vert<R$ implies 
 $ \frac 9 {10} (1+\vert x\vert) \leq 1+\vert y\vert\leq \frac{11}{10} (1+\vert x\vert)$, thus we can control the left
side in (\ref{muck}) for $w_\delta$  by $( \frac {11}{9})^{\frac\delta p}$. 
 
 Finally, for $R>1$ and $\vert x\vert\leq 10 R$, we write
 \begin{equation*} \begin{split}
  & \left(\frac 1{\vert B(x,R)\vert}  \!\int_{B(x,R)}\! \!\!\!\! w(y)\, dy\right)^{\frac 1 p} \left( \frac 1{\vert B(x,R)\vert}\! \int_{B(0,R)} \!\frac{dy}{w(y)^{\frac 1{p-1}}}\right)^{1-\frac 1 p} 
  \\ \leq &   \left(\frac 1{\vert B(0, R)\vert}  \!\int_{B(x,11\, R)}\! \!\!\!\!w(y)\, dy\right)^{\frac 1 p}\!\! \left( \frac 1{\vert B(0,R)\vert}\! \int_{B(0,11 \,R)} \!\frac{dy}{w(y)^{\frac 1{p-1}}}\right)^{1-\frac 1 p} 
  \\ = & \left(\frac 1 {R^3} \int_0^{11\, R} r^2 \frac{dr}{(1+r)^\delta}\right)^{\frac 1 p} \left(\frac 1 {R^3} \int_0^{11\, R} r^2(1+r)^{\frac \delta {p-1}}\, dr\right)^{1-\frac 1 p}
    \\ \leq &  c_{\delta, p} \left(\frac 1 {R^3} \int_0^{11\, R} r^2 \frac{dr}{ r^\delta}\right)^{\frac 1 p} \left( \left(\frac 1 {R^3} \int_0^{11\, R} r^2 dr \right)^{1-\frac 1 p} + \left( \frac 1 {R^3} \int_0^{11\, R} r^{2+\frac \delta {p-1}}\, dr\right)^{1-\frac 1 p} \right)
    \\ = & c_{\delta, p} \frac{11^3} { (3-\delta)^{\frac 1 p} } \left( \frac{ (11R)^{ -\frac{\delta}{p} } }{  3^{1-\frac 1 p}} +  \frac{ 1  }{ (3+\frac \delta{p-1})^{1-\frac 1 p}} \right).
  \end{split}\end{equation*}
  The lemma is proved. \Endproof
  
\begin{lemma}
  \label{LemRieszt}
    If $0<\delta<3$ and $1<p<+\infty$, then the Riesz transforms $R_i$ and the Hardy--Littlewood maximal function operator are bounded on $L^p_{w_\delta}=L^p(w_\delta(x)\, dx)$ :
 $$ \|R_jf\|_{L^p_{w_\delta}}\leq C_{p,\delta} \|f\|_{L^p_{w_\delta}} \text{ and } \|\mathcal{M}_f\|_{L^p_{w_\delta}}\leq C_{p,\delta} \|f\|_{L^p_{w_\delta}}.
 $$
  \end{lemma}
  
  \Proof{}
 The boundedness of the Riesz transforms or of the Hardy--Littlewwod maximal function on $L^p(w_\gamma\, dx)$ are   basic properties of the Muckenhoupt class $\mathcal{A}_p$ \cite{Gr09}.
  \Endproof{}\\
  
  \noindent
  We will use strategically the next corollary, which is specially useful
  to obtain discretely self-similar solutions.
  
  \begin{corollary}[Non-increasing kernels] Let $\theta\in L^1(\mathbb{R}^3)$ be a non-negative radial function which is radially non-increasing. Then, if $0<\delta<3$ and $1<p<+\infty$, we have, for $f\in L^p_{w_\delta}$, the inequality
  $$ \| \theta*f\|_{L^p_{w_\delta}} \leq C_{p,\delta} \|f\|_{L^p_{w_\delta}} \|\theta\|_1.$$
  \end{corollary}
  
  \Proof We have the well-known inequality for radial non-increasing kernels \cite{Gr08}
  $$ \vert \theta*f(x)\vert\leq \|\theta\|_1 \mathcal{M}_f(x)$$ so that we may conclude with Lemma \ref{LemRieszt}. \Endproof\\

  \noindent
  We illustrate the utility of Lemma \ref{LemRieszt}  with the following corollaries:
  
  \begin{corollary}
  \label{CorRieszt}
  Let $0<\gamma < \frac{5}{2} $ and $0<T<+\infty$. Let   $\mathbb{F}$ be  a tensor $\mathbb{F}(t,x)=\left(F_{i,j}(t,x)\right)_{1\leq i,j\leq 3}$ such that $\mathbb{F}\in L^2((0,T), L^2_{w_\gamma})$.  Let $\vb$ be a time-dependent divergence free vector-field ($\vN\cdot\vb=0$) such that $\vb\in L^3((0,T),L^3_{w_{3\gamma/2}})$.
  
 Let $\vu$ be a solution of the following advection-diffusion problem 
 \begin{equation}  \label{advection difusion free} \left\{ \begin{matrix} \partial_t \vu= \Delta \vu  -(\vb\cdot \vN)\vu- \vN q +\vN\cdot \mathbb{F} \cr \cr \vN\cdot \vu=0,    \phantom{space space} 
 \end{matrix}\right.\end{equation}
be such that :
  $\vu$ belongs to $L^\infty((0,T), L^2_{w_\gamma})$ and $\vN\vu$ belongs to $L^2((0,T),L^2_{w_\gamma})$, and the pressure $q$ belongs to $\mathcal{D}'( (0,T) \times \mathbb{R} ^3 )$. \\
  
  Then, the gradient of the pressure $\vN q$ is necessarily
  related to $\vu$, $\vb$  and $\mathbb{F}$ through the Riesz transforms $R_i =\frac{\partial_i}{\sqrt{-\Delta}}$ by the formula
 $$ \vN q= \vN \left( \sum_{i=1}^3\sum_{j=1}^3 R_iR_j(b_iu_j-F_{i,j}) \right) $$
 and  $\sum_{i=1}^3\sum_{j=1}^3 R_iR_j(b_iu_j)$  belongs to $L^{3}((0,T),L^{6/5}_{w_{\frac {6\gamma}5}})$ and  $\sum_{i=1}^3\sum_{j=1}^3 R_iR_j F_{i,j}  $ belongs to $L^{2}((0,T),L^{2}_{w_\gamma})$.
  
  \end{corollary}
  
  \Proof{}
  We define 
  \begin{equation*}
      p =  \left( \sum_{i=1}^3\sum_{j=1}^3 R_iR_j(b_iu_j-F_{i,j}) \right).
  \end{equation*}
  As $ 0 < \gamma < \frac{5}{2} $ we can use Lemma \ref{LemRieszt} and \eqref{productbu} to obtain
  $\sum_{i=1}^3\sum_{j=1}^3 R_iR_j(b_i u_j)$  belongs to $L^{3}((0,T),L^{6/5}_{w_{\frac {6\gamma}5}})$ and  $\sum_{i=1}^3\sum_{j=1}^3 R_iR_j F_{i,j}  $ belongs to $L^{2}((0,T),L^{2}_{w_\gamma})$.

  Taking the divergence in \eqref{advection difusion free}, we obtain $ \Delta (q-p)=0 $. We take a test function $\alpha \in \mathcal{D}(\mathbb{R})$ such that $\alpha (t)= 0$  for all $|t| \geq \varepsilon$, and a test function $\beta\in\mathcal{D}(\mathbb{R}^3)$;  then the distribution $\vN q *( \alpha\otimes\beta) $ is well defined on $(\varepsilon, T-\varepsilon) \times \mathbb{R}^3$.
  
  We fix  $t \in (\varepsilon, T-\varepsilon)$ and define $$A_{\alpha,\beta,t}=( \vN q * (\alpha\otimes\beta)- \vN p * (\alpha\otimes\beta)) (t,.). $$  We have 
 \begin{equation}\label{convol} \begin{split} A_{\alpha,\beta,t}=& ( \vu *(-\partial_t\alpha\otimes\beta+\alpha\otimes\Delta\beta)  +(-\vu\otimes\vb +\mathbb{F}) \cdot (\alpha\otimes\vN\beta))(t,.)\\&  -  ( p * (\alpha\otimes\vN\beta)) (t,.) .\end{split} \end{equation} 
 Convolution with a function in $\mathcal{D}(\mathbb{R}^3)$ is a bounded operator on $L^2_{w_\gamma}$ and on $L^{6/5}_{w_{6\gamma/5}}$ (as, for $\varphi\in \mathcal{D}(\mathbb{R}^3)$ we have $\vert f*\varphi\vert\leq C_\varphi \mathcal{M}_f$). Thus, we may conclude from (\ref{convol}) that $A_{\alpha,\beta,t}\in L^2_{w_\gamma}+ L^{6/5}_{w_{6\gamma/5}}$. If $\max \{ \gamma, \frac{\gamma  + 2 }{2} \}  <\delta<5/2$ , we have $A_{\alpha,\beta,t}\in L^{6/5}_{w_{6\delta/5}}$.
  
  In particular,  $A_{\alpha,\beta,t}$ is a tempered  distribution.    As we have $$\Delta  A_{\alpha,\beta,t}=(\alpha\otimes\beta)*(\Delta(q-p))(t,.)=0,$$ we find that  $A_{\alpha,\beta,t}$ is a polynomial. 
   We remark that for all $1<r<+\infty$ and $0< \delta < 3$, $L^r_{w_\delta}$ does not contain non-trivial   polynomials. 
  Thus, $ A_{\alpha,\beta,t}= 0$. We then use an approximation of identity $\frac 1{\epsilon^4} \alpha	(\frac t\epsilon)\beta(\frac x\epsilon)$ and conclude that $\vN(q-p)=0$.  \Endproof{}

  Actually, we can answer a question posed by Bradshaw and Tsai in \cite{BT19} about the nature of the pressure for self-similar solutions of the Navier--Stokes equations. In effect, we have the next corollary: 
  
 $\ $
  
  \begin{corollary}
  \label{pressurenature2}
  Let $1<\gamma < \frac{5}{2}$ and $0<T<+\infty$.  Let   $\mathbb{F}$ be  a  tensor $\mathbb{F}(t,x)=\left(F_{i,j}(t,x)\right)_{1\leq i,j\leq 3}$ such that $\mathbb{F}\in L^2((0,T), L^2_{w_\gamma})$.
  
 Let $\vu$ be a solution of the following problem 
 \begin{equation*}  \label{advection difusion free2} \left\{ \begin{matrix} \partial_t \vu= \Delta \vu  -(\vu \cdot \vN)\vu- \vN p +\vN\cdot \mathbb{F} \cr \cr \vN\cdot \vu=0,    \phantom{space space} 
 \end{matrix}\right.\end{equation*}
be such that :
  $\vu$ belongs to $L^\infty([0,+\infty ), L^2)_{loc}$ and $\vN \vu$ belongs to $L^2([0,+\infty ),L^2)_{loc}$, and the pressure $q$ is in $\mathcal{D}'( (0,T) \times \mathbb{R} ^3 )$. \\
  
  We suppose that there exists  $\lambda>1$ such that $ \lambda^2 \mathbb{F}(\lambda^2 t,\lambda x)=\mathbb{F}(t,x) $ and  $\lambda \vu(\lambda^2 t,\lambda x)=\vu(t,x)$.
  Then, the gradient of the pressure $\vN q$ is necessarily
  related to $\vu$  and $\mathbb{F}$ through the Riesz transforms $R_i =\frac{\partial_i}{\sqrt{-\Delta}}$ by the formula
 $$ \vN q= \vN \left( \sum_{i=1}^3\sum_{j=1}^3 R_iR_j(u_iu_j-F_{i,j}) \right) $$
 and  $\sum_{i=1}^3\sum_{j=1}^3 R_iR_j(u_iu_j)$  belongs to $L^{4}((0,T),L^{6/5}_{w_{\frac {6\gamma}5}})$ and  $\sum_{i=1}^3\sum_{j=1}^3 R_iR_j F_{i,j}  $ belongs to $L^{2}((0,T),L^{2}_{w_\gamma})$.
  \end{corollary}
  
  \Proof{}
  We shall  use Corollary \ref{CorRieszt}, and thus we need to show that $\vu$ belongs to $L^\infty((0,T), L^2_{w_\gamma} \cap L^3((0,T), L^3_{3\gamma/2}))$ and $\vN\vu$ belongs to $L^2((0,T),L^2_{w_\gamma})$. In fact, 
  \begin{equation*}
      \| u \|_{L^\infty ((0,T),L^2_{w_\gamma})} \leq \sup_{0 \leq t \leq T} \int_{|x|<1} |\vu (t,x)|^2 \, dx + c \sup_{0 \leq t \leq T} \sum_{k\in \mathbb{N}}  \int_{\lambda^{k-1}<|x|<\lambda^k} \frac{|\vu (t,x)|^2}{\lambda^{\gamma k } } \, dx
  \end{equation*}
  and
  \begin{align*}
     \sup_{0 \leq t \leq T} \sum_{k\geq 1}  \int_{\lambda^{k-1}<|x|<\lambda^k} \frac{|\vu (t,x)|^2}{\lambda^{\gamma k } } \, dx & \leq \sup_{0 \leq t \leq T} \sum_{k\in \mathbb{N}} \lambda^{ (1-\gamma) k } \int_{\lambda^{-1}<|x|< 1} |\vu (\frac{t}{\lambda^{2k}},x)|^2 \, dx \\
     & \leq  c \sup_{0 \leq t \leq T} \int_{\lambda^{-1}<|x|< 1} |\vu (t,x)|^2 \, dx < +\infty.
  \end{align*}
  For $\vN\vu$, we compute
  for $k\in\mathbb{N}$, 
 $$\int_0^T \int_{\lambda^{k-1}<\vert x\vert<\lambda^{k}} \vert \vN \vu(t,x)\vert^2\, dt\, dx=\lambda^{k}\int_0^{\frac T{\lambda^{2k}}}\int_{\frac 1 \lambda <\vert x\vert<1} \vert \vN \vu(t,x)\vert^2\, dx\, dt.
 $$ We may conclude that $\vN\vu$ belongs to $L^2((0,T),L^2_{w_\gamma})$, since for $\gamma>1$ we have $\sum_{k\in\mathbb{N}} \lambda^{ (1-\gamma) k }<+\infty$.
 
 Now, we use the Sobolev embeddings described in next Lemma (Lemma \ref{sobol}) to get that $\vu$ belongs to $L^2((0,T),L^6_{w_{3\gamma}})$, and thus (by interpolation with $L^\infty((0,T),L^2_{w_\gamma}))$ to $L^4((0,T),L^3_{w_{3\gamma/2}})$.
 
 \noindent
 In particular, $\sum_{i=1}^3\sum_{j=1}^3 R_iR_j(u_iu_j)$  belongs to $L^{4}((0,T),L^{6/5}_{w_{\frac {6\gamma}5}})$, since we have 
 \begin{equation*}
        \label{productbu}
      \|  ( \vu \otimes \vu ) w_\gamma   \|_{ L^{6/5}} \leq \| \sqrt{w_\gamma} \vu \|_{ L^{2} } \| \sqrt{w_\gamma} \vu \|_{ L^{3} } \leq \| \sqrt{w_\gamma} \vu \|^{\frac{3}{2}}_{ L^{2} } \| \sqrt{w_\gamma} \vu \|^{\frac{1}{2}}_{ L^{6} .}
  \end{equation*} 
  \Endproof{}

  \begin{lemma}[Sobolev embeddings] \label{sobol} Let $\delta>0$. If $f\in L^2_{w_\delta}$ and $\vN f\in L^2_{w_\delta}$ then $f\in L^6_{w_{3\delta}}$ and
  $$ \|f\|_{L^6_{w_{3\delta}}}\leq C_\delta (\|f\|_{L^2_{w_\delta}}+ \|\vN f\|_{L^2_{w_\delta}}).$$
  \end{lemma}
  
  \Proof Since both $f$ and $w_{\delta/2}$ are locally in $H^1$, we write
  $$ \partial_i(f w_{\delta/2})=w_{\delta/2} \partial_i f+ f \partial_i(w_{\delta/2})= w_{\delta/2} \partial_i f -\frac\delta 2 \frac {x_i}{\vert x\vert} w_{\delta/2} f$$
  and thus
  $$ \| w_{\delta/2}f\|_2^2+ \|\vN(w_{\delta/2}f)\|_2^2\leq  (1+\frac{\delta^2}2) \| w_{\delta/2}f\|_2^2 + 2 \| w_{\delta/2}\vN f\|_2^2.$$
  Thus, $w_{\delta/2} f$ belongs to $L^6$ (since $H^1\subset L^6$), or equivalently $f\in L^6_{w_{3\delta}}$.\Endproof
 
\section{A priori estimates for the advection-diffusion problem.}

\subsection{Proof of Theorem \ref{estimates}.}
Let $0<t_0<t_1<T$. We take a function $
\alpha\in \mathcal{C}^\infty(\mathbb{R})$ which is non-decreasing, with $\alpha(t)$ equal to $0$ for $t<1/2$ and equal to $1$ for $t>1$.  For $0<\eta< \min(\frac {t_0}2,T-t_1) $, we define $$\alpha_{\eta,t_0,t_1}(t)=\alpha( \frac{t-t_0}\eta)-\alpha(\frac{t-t_1}\eta) .$$ We take as well a non-negative  function $\phi\in\mathcal{D}(\mathbb{R}^3)$ which is equal to $1$ for $\vert x\vert\leq 1$ and to $0$ for $\vert x\vert\geq 2$. For $R>0$, we define $\phi_R(x)=\phi(\frac x R)$. Finally, we define, for $\epsilon>0$, $w_{\gamma,\epsilon}= \frac 1{(1+\sqrt{\epsilon^2+\vert x\vert^2})^\delta}$. We have $\alpha_{\eta,t_0,t_1}(t)\phi_R(x) w_{\gamma,\epsilon}(x)\in\mathcal{D}((0,T)\times\mathbb{R}^3)$ and $\alpha_{\eta,t_0,t_1}(t)\phi_R(x) w_{\gamma,\epsilon}(x) \geq 0$. Thus, using the local energy balance (\ref{energloc}) and the fact that $\mu\geq 0$, we find
\begin{equation*} \begin{split}
 -\iint \frac{\vert\vu\vert^2}2& \partial_t\alpha_{\eta,t_0,t_1} \phi_R w_{\gamma,\epsilon}\, dx\, ds
\\ \leq&-\sum_{i=1}^3 \iint \partial_i\vu\cdot \vu\,  \alpha_{\eta,t_0,t_1}  (w_{\gamma,\epsilon}\partial_i \phi_R+\phi_R \partial_iw_{\gamma,\epsilon})\, dx\, ds
\\&- \iint \vert\vN\vu\vert^2\, \,  \alpha_{\eta,t_0,t_1} \phi_R w_{\gamma,\epsilon} dx\, ds
\\&  + \sum_{i=1}^3 \iint  \frac{\vert\vu\vert^2}2 b_i  \alpha_{\eta,t_0,t_1}  (w_{\gamma,\epsilon}\partial_i \phi_R+\phi_R \partial_iw_{\gamma,\epsilon})\, dx\, ds
\\& +\sum_{i=1}^3 \iint     \alpha_{\eta,t_0,t_1} pu_i (w_{\gamma,\epsilon}\partial_i \phi_R+\phi_R \partial_iw_{\gamma,\epsilon})\, dx\, ds
\\&- \sum_{i=1}^3\sum_{j=1}^3 \iint  F_{i,j} u_j  \alpha_{\eta,t_0,t_1}  (w_{\gamma,\epsilon}\partial_i \phi_R+\phi_R \partial_iw_{\gamma,\epsilon})\, dx\, ds
\\&  - \sum_{i=1}^3 \sum_{j=1}^3\iint   F_{i,j}\partial_iu_j\  \alpha_{\eta,t_0,t_1} \phi_R w_{\gamma,\epsilon}\, dx\, ds.
\end{split}\end{equation*}
We remark that, independently from $R>1$ and $\epsilon>0$, we have (for $0<\gamma\leq 2$)
$$ \vert  w_{\gamma,\epsilon}\partial_i \phi_R\vert +\vert \phi_R \partial_iw_{\gamma,\epsilon}\vert \leq C_\gamma \frac{w_\gamma(x)}{1+\vert x\vert}\leq C_\gamma w_{3\gamma/2}(x).
$$
Moreover, we know that $\vu$ belongs to $L^\infty((0,T),L^2_{w_\gamma})\cap L^2((0,T), L^6_{w_{3\gamma}})$ hence to $L^4((0,T),L^3_{w_{3\gamma/2}})$. Since $T<+\infty$, we have as well $\vu\in L^3((0,T), L^3_{w_{3\gamma/2}})$. (This is the same type of integrability as required for $\vb$).
 Moreover, we have $p u_i\in L^1_{w_{3\gamma/2}}$ since $w_\gamma p\in L^2 ((0,T), L^{6/5}+L^2)$ and $w_{\gamma/2} \vu\in L^2((0,T), L^2\cap L^6)$. All those remarks will allow us to use dominated convergence.

 We first let $\eta$ go to $0$. We find that 
 
\begin{equation*} \begin{split}
- \lim_{\eta\rightarrow 0}  \iint \frac{\vert\vu\vert^2}2& \partial_t\alpha_{\eta,t_0,t_1}   \phi_R w_{\gamma,\epsilon}\, dx\, ds
 \\ \leq&-\sum_{i=1}^3 \int_{t_0}^{t_1}\int   \partial_i\vu\cdot \vu\,  (w_{\gamma,\epsilon}\partial_i \phi_R+\phi_R \partial_iw_{\gamma,\epsilon})\, dx\, ds
\\&- \int_{t_0}^{t_1} \int  \vert\vN\vu\vert^2\, \,  \phi_R w_{\gamma,\epsilon} dx\, ds
\\&  + \sum_{i=1}^3 \int_{t_0}^{t_1}\int  \frac{\vert\vu\vert^2}2 b_i     (w_{\gamma,\epsilon}\partial_i \phi_R+\phi_R \partial_iw_{\gamma,\epsilon})\, dx\, ds
\\& +\sum_{i=1}^3 \int_{t_0}^{t_1}   \int     pu_i (w_{\gamma,\epsilon}\partial_i \phi_R+\phi_R \partial_iw_{\gamma,\epsilon})\, dx\, ds
\\&- \sum_{i=1}^3\sum_{j=1}^3 \int_{t_0}^{t_1}\int    F_{i,j} u_j    (w_{\gamma,\epsilon}\partial_i \phi_R+\phi_R \partial_iw_{\gamma,\epsilon})\, dx\, ds
\\&  - \sum_{i=1}^3 \sum_{j=1}^3\int_{t_0}^{t_1} \int    F_{i,j}\partial_iu_j\    \phi_R w_{\gamma,\epsilon}\, dx\, ds.
\end{split}\end{equation*}
 
 Let us define
 $$ A_{R,\epsilon}(t)=\int \vert \vu(t,x)\vert^2   \phi_R(x)  w_{\gamma,\epsilon}(x)\, dx.$$
As we have
$$-\iint \frac{\vert\vu\vert^2}2 \partial_t\alpha_{\eta,t_0,t_1}   \phi_R w_{\gamma,\epsilon} \, dx\, ds=-\frac 1 2\int \partial_t\alpha_{\eta,t_0,t_1}
A_{R,\epsilon}(s) \, ds$$
we find that, when $t_0$ and $t_1$ are Lebesgue points of the measurable function $A_{R,\epsilon}$ $$ \lim_{\eta\rightarrow 0}  -\iint \frac{\vert\vu\vert^2}2 \partial_t\alpha_{\eta,t_0,t_1}  \phi_R w_{\gamma,\epsilon} \, dx\, ds=\frac 1 2 (  A_{R,\epsilon}(t_1)- A_{R,\epsilon}(t_0))   .$$ Then, by continuity, we can let $t_0$ go to $0$ and thus  replace $t_0$ by $0$ in the inequality. Moreover, if we let $t_1$ go to $t$, then by weak continuity,
we find that $ A_{R,\epsilon}(t)\leq \lim_{t_1\rightarrow t }  A_{R,\epsilon}(t_1)$, so that we may as well replace $t_1$ by $t\in (0,T)$. Thus we find that for every $t\in (0,T)$, we have 

 \begin{equation}\label{localineq} \begin{split}
\int \frac{\vert \vu(t,x)\vert^2}2 \phi_R w_{\gamma,\epsilon}\, dx     
 \\ \leq& \int \frac{\vert \vu_0(x)\vert^2}2 \phi_R w_{\gamma,\epsilon}\, dx     
 \\ &-\sum_{i=1}^3 \int_0^t \int   \partial_i\vu\cdot \vu\,  (w_{\gamma,\epsilon}\partial_i \phi_R+\phi_R \partial_iw_{\gamma,\epsilon})\, dx\, ds
\\&- \int_0^t  \int  \vert\vN\vu\vert^2\, \,  \phi_R w_{\gamma,\epsilon} dx\, ds
\\&  + \sum_{i=1}^3 \int_0^t \int  \frac{\vert\vu\vert^2}2 b_i     (w_{\gamma,\epsilon}\partial_i \phi_R+\phi_R \partial_iw_{\gamma,\epsilon})\, dx\, ds
\\& +\sum_{i=1}^3 \int_0^t    \int     pu_i (w_{\gamma,\epsilon}\partial_i \phi_R+\phi_R \partial_iw_{\gamma,\epsilon})\, dx\, ds
\\&- \sum_{i=1}^3\sum_{j=1}^3 \int_0^t \int    F_{i,j} u_j    (w_{\gamma,\epsilon}\partial_i \phi_R+\phi_R \partial_iw_{\gamma,\epsilon})\, dx\, ds
\\&  - \sum_{i=1}^3 \sum_{j=1}^3\int_0^t  \int    F_{i,j}\partial_iu_j\    \phi_R w_{\gamma,\epsilon}\, dx\, ds.
\end{split}\end{equation}

Thus, letting $R$ go to $+\infty$ and then $\epsilon$ go to $0$, we find by dominated convergence that,
  for every $t\in (0,T)$, we have
 \begin{equation*}\begin{split}
 \|\vu(t,. )\|_{L^2_{w_\gamma}}^2 &+2\int_0^t \|\vN \vu(s,.)\|_{L^2_{w_\gamma}}^2\, ds
 \\  \leq& \|\vu_0 \|_{L^2_{w_\gamma}}^2 -\int_0^t \int   \vN\vert\vu\vert^2\cdot\vN w_\gamma\, dx\, ds+\int_0^t\int (\vert\vu\vert^2 \vb+2p\vu) \cdot \vN(w_\gamma)\, dx\, ds\\ &\phantom{space space} -2\sum_{i=1}^3\sum_{j=1}^3 \int_0^t\int F_{i,j} (\partial_i u_j) w_\gamma +F_{i,j} u_i \partial_j(w_\gamma)\, dx\, ds.
 \end{split}\end{equation*}
 Now we write
 \begin{equation*}\begin{split} \left\vert \int_0^t\int \vN\vert \vu\vert^2\cdot \vN w_\gamma\, ds\, ds\right\vert\leq&  2\gamma \int_0^t\int  \vert\vu\vert \vert\vN\vu\vert \, w_\gamma\, dx\, ds \\ \leq& \frac 1 4 \int_0^t  \|\vN\vu\|_{L^2_{w_\gamma}}^2\, ds+4\gamma^2 \int_0^t \|\vu\|_{L^2_{w_\gamma}}^2\, ds .
 \end{split}\end{equation*}
Writing $$ p_1=\sum_{i=1}^3\sum_{j=1}^3 R_iR_j(b_iu_j) \text{ and }   p_2=-\sum_{i=1}^3\sum_{j=1}^3 R_iR_j(F_{i,j})$$
and using the fact that $w_{6\gamma/5}\in  \mathcal{A}_{6/5}$ and $w_\gamma\in\mathcal{A}_2$, we get

 \begin{equation*}\begin{split} \left\vert \int_0^t\int  (\vert\vu\vert^2 \vb+2p_1\vu) \cdot \vN(w_\gamma)\, dx\, ds\right\vert\leq&   \gamma \int_0^t\int   (\vert\vu\vert ^2 \vert  \vb\vert+2\vert p_1\vert \,\vert\vu\vert) \, w_\gamma^{3/2}\, dx\, ds \\ \leq   \gamma \int_0^t  &\|w_\gamma^{1/2} \vu\|_6 (\|  w_\gamma \vert \vb\vert \vert\vu\vert  \| _{6/5}+\|w_\gamma p_1\|_{6/5})  ds \\ \leq 
  C_\gamma \int_0^t & \|w_\gamma^{1/2} \vu\|_6  \|  w_\gamma \vert \vb\vert \vert\vu\vert \|_{6/5}  \, ds \\ \leq
  C_\gamma \int_0^t & \|w_\gamma^{1/2} \vu\|_6  \|  w_\gamma^{1/2} \vb\|_3    \|w_\gamma^{1/2}\vu\|_2    \, ds \\\leq 
  C_\gamma'  \int_0^t  & (\|\vN\vu\|_{L^2_{w_\gamma}}+\|\vu\|_{L^2_{w_\gamma}})\  \|\vb\|_{L^3_{w_{3\gamma/2}}}\|\vu\|_{L^2_{w_\gamma}}    \, ds \\ \leq \frac 1 4 \int_0^t  \|\vN\vu\|_{L^2_{w_\gamma}}^2\, ds &+ C''_\gamma \int_0^t \|\vu\|_{L^2_{w_\gamma}}^2 ( \|\vb\|_{L^3_{w_{3\gamma/2}}} + \|\vb\|_{L^3_{w_{3\gamma/2}}}^2)\, ds
 \end{split}\end{equation*}
and  
 \begin{equation*}\begin{split} \left\vert \int_0^t\int   2p_2\vu \cdot \vN(w_\gamma)\, dx\, ds\right\vert\leq&  2  \gamma \int_0^t\int  \vert p_2\vert\, \vert \vu\vert \, w_\gamma \, dx\, ds
 \\ \leq& \gamma \int_0^t \|\vu\|_{L^2_{w_\gamma}}^2+\|p_2\|_{L^2_{w_\gamma}}^2\, ds \\ \leq    & C_\gamma \int_0^t \|\vu\|_{L^2_{w_\gamma}}^2+\|\mathbb{F}\|_{L^2_{w_\gamma}}^2\, ds. 
 \end{split}\end{equation*}
 Finally, we have
  \begin{equation*}\begin{split} \left\vert2\sum_{i=1}^3\sum_{j=1}^3 \int_0^t\int F_{i,j} (\partial_i u_j) w_\gamma +F_{i,j} u_i \partial_j(w_\gamma)\, dx\, ds\right\vert\leq&  2    \int_0^t\int  \vert  F\vert\, (\vert \vN\vu\vert+\gamma \vert \vu\vert) \, w_\gamma \, dx\, ds \\ \leq    \frac 1 4 \int_0^t  \|\vN\vu\|_{L^2_{w_\gamma}}^2\, ds &+ C_\gamma \int_0^t \|\vu\|_{L^2_{w_\gamma}}^2+\|\mathbb{F}\|_{L^2_{w_\gamma}}^2\, ds. 
 \end{split}\end{equation*}
We have obtained
 \begin{equation}\label{gronwb} \begin{split}  \|\vu(t,. )\|_{L^2_{w_\gamma}}^2 &+  \int_0^t  \|\vN\vu\|_{L^2_{w_\gamma}}^2\, ds  \\ \leq& \|\vu_0 \|_{L^2_{w_\gamma}}^2  +C_\gamma \int_0^t \|  \mathbb{F}(s,.)\|^2_{L^2_{w_\gamma}}\, ds + C_\gamma  \int_0^t  (1+  \|\vb(s,. )\|_{L^3_{w_{3\gamma/2}}}^2 )  \|\vu(s,. )\|_{L^2_{w_\gamma}}^2  \, ds\end{split}\end{equation} 
 and Theorem \ref{estimates} is proven.\Endproof

 \subsection{Passive transportation.}
 From inequality (\ref{gronwb}), we have the following direct consequence :
 \begin{corollary}\label{passive} Under the assumptions of Theorem \ref{estimates}, we have
 $$ \sup_{0<t<T} \|\vu\|_{L^2_{w_\gamma}} \leq (\|\vu_0\|_{L^2_{w_\gamma}}+ C_\gamma \|\mathbb{F}\|_{L^2((0,T), L^2_{w_\gamma})}) \ e^{C_\gamma (T+ T^{1/3} \|\vb\|_{L^3((0,T), L^3_{w_{3\gamma /2}})}^2)}$$
 and
 $$  \|\vN\vu\|_{L^2((0,T),L^2_{w_\gamma)}} \leq (\|\vu_0\|_{L^2_{w_\gamma}}+ C_\gamma \|\mathbb{F}\|_{L^2((0,T), L^2_{w_\gamma})}) \ e^{C_\gamma (T+ T^{1/3} \|\vb\|_{L^3((0,T), L^3_{w_{3\gamma /2}})}^2)}$$ where the constant $C_\gamma$ depends only on $\gamma$.
 \end{corollary}
 Another direct consequence is the following uniqueness result for the advection-diffusion problem with a (locally in time), bounded $\vb$ :
 
  \begin{corollary}\label{unique}. 
  Let $0<\gamma\leq 2$. Let $0<T<+\infty$. Let  $\vu_{0}$ be a divergence-free vector field such that $\vu_0\in L^2_{w_\gamma}(\mathbb{R}^3)$ and   $\mathbb{F}$ be  a tensor $\mathbb{F}(t,x)=\left(F_{i,j}(t,x)\right)_{1\leq i,j\leq 3}$ such that $\mathbb{F}\in L^2((0,T), L^2_{w_\gamma})$.  Let $\vb$ be a time-dependent divergence free vector-field ($\vN\cdot\vb=0$) such that $\vb\in L^3((0,T),L^3_{w_{3\gamma/2}})$. Assume moreover that $\vb$ belongs   to $L^2_t L^\infty_x(K)$ for every compact subset $K$ of  $(0,T)\times\mathbb{R}^3 $. 
  
 Let $(\vu_1, p_1)$ and $(\vu_2,p_2)$  be two solutions of the following advection-diffusion problem 
 \begin{equation*}  (AD) \left\{ \begin{matrix} \partial_t \vu= \Delta \vu  -(\vb\cdot \vN)\vu- \vN p +\vN\cdot \mathbb{F} \cr \cr \vN\cdot \vu=0,    \phantom{space space} \vu(0,.)=\vu_0
 \end{matrix}\right.\end{equation*}
be such that, for $k=1$ and $k=2$,  :
 \begin{itemize} 
 \item[$\bullet$]  $\vu_k$ belongs to $L^\infty((0,T), L^2_{w_\gamma})$ and $\vN\vu_k$ belongs to $L^2((0,T),L^2_{w_\gamma})$
 \item[$\bullet$] the pressure $p_k$  is related to $\vu_k$, $\vb$  and $\mathbb{F}$ through the Riesz transforms $R_i =\frac{\partial_i}{\sqrt{-\Delta}}$ by the formula
 $$ p_k=\sum_{i=1}^3\sum_{j=1}^3 R_iR_j(b_iu_{k,j}-F_{i,j})$$
 \item[$\bullet$] the map $t\in [0,T)\mapsto \vu_k(t,.)$ is weakly continuous from $[0,T)$ to $L^2_{w_\gamma}$, and is strongly continuous at $t=0$ :
 $$ \lim_{t\rightarrow 0} \|\vu_k(t,.)-\vu_0\|_{L^2_{w_\gamma}}=0.$$
 \end{itemize}
 
 Then $\vu_1=\vu_2$.
    \end{corollary}
    
    \Proof Let $\vv=\vu_1-\vu_2$ and $q=p_1-p_2$. Then we have
     \begin{equation*}   \left\{ \begin{matrix} \partial_t \vv= \Delta \vv  -(\vb\cdot \vN)\vv- \vN q   \cr \cr \vN\cdot \vv=0,    \phantom{space space} \vv(0,.)=0
 \end{matrix}\right.\end{equation*}
Moreover on every compact subset $K$ of $(0,T)\times\mathbb{R}^3$, $\vb\otimes\vv$ is in $L^2_t L^2_x$, while it belongs globally to $L^{3}_t L^{6/5}_{w_{6\gamma/5}}$. Writing, for $\varphi, \psi\in\mathcal{D}((0,T)\times\mathbb{R}^3)$ such that $\psi=1$ on the neigborhood of the support of $\varphi$,
$$ \varphi q=q_1+q_2=\varphi \sum_{i=1}^3\sum_{j=1}^3 R_iR_j(\psi b_iv_j)+\varphi \sum_{i=1}^3\sum_{j=1}^3 R_iR_j((1-\psi) b_i v_j)$$ we find that $\|q_1\|_{L^2L^2}\leq C_{\varphi,\psi} \|\psi \vb\otimes\vv\|_{L^2 L^2}$ and
$$ \|q_2\|_{L^3 L^\infty} \leq C_{\varphi,\psi} \|\vb\otimes \vv\|_{L^3 L^{6/5}_{w_{6\gamma/5}}}
$$
with
$$ C_{\varphi,\psi}\leq C \|\varphi\|_\infty \|1-\psi\|_\infty  \sup_{x\in {\rm Supp}\, \varphi}  \left( \int_{y\in {\rm Supp }\, (1-\psi)}  \left( \frac { (1+\vert y\vert)^\gamma} {\vert x-y\vert^3}\right)^6 \right)^{1/6}<+\infty.$$
Thus,  we may take the scalar product of  $\partial_t \vv$ with $\vv$ and find  that
 \begin{equation*}\label{} \partial_t(\frac {\vert\vv\vert^2}2)=\Delta(\frac {\vert\vv\vert^2}2)-\vert\vN \vv\vert^2- \vN\cdot\left( \frac{\vert\vv\vert^2}2\vb\right)-\vN\cdot(q\vv) . \end{equation*}
 Thus we are under the assumptions of Theorem \ref{estimates} and we may use Corollary \ref{passive} to find that $\vv=0$. \Endproof

\subsection{Active transportation.}
We begin with the following lemma :
\begin{lemma}\label{gronwallnl} Let $\alpha$ be a non-negative bounded measurable  function on $[0,T)$ such that, for two constants $A,B\geq 0$, we have
$$ \alpha(t)\leq A + B\int_0^t \alpha(s)+\alpha(s)^3\, ds.$$ If $T_0>0$ and $T_1=\min(T,T_0, \frac 1{4B (A+BT_0)^2})$, we have, for every $t\in [0,T_1]$,  $\alpha(t)\leq \sqrt{ 2} (A+BT_0)$.
\end{lemma}
\Proof We write $\alpha\leq 1+\alpha^3$. We define 
$$\Phi(t)= A+B T_0 + B \int_0^t \alpha^3\, ds \text{ and } \Psi(t)=A+BT_0+ B \int_0^t  \Phi^3(s)\, ds.$$ 
We have, for $t\in [0,T_1]$,  $\alpha\leq \Phi\leq \Psi$. Since $\Psi$ is $\mathcal{C}^1$, we may write
$$ \Psi'(t)= B \Phi(t)^3 \leq B \Psi(t)^3$$
and thus
$$ \frac 1 {\Psi(0)^2}-\frac 1{\Psi(t)^2}\leq 2Bt.$$ We thus find
$$ \Psi(t)^2 \leq \frac{\Psi(0)^2}{1-2B \Psi(0)^2 t}\leq 2 \Psi(0)^2.$$
The lemma is proven.\Endproof

\begin{corollary}\label{active}  Assume that $\vu_0$, $\vu$, $p$, $\mathbb{F}$ and $\vb$ satisfy  assumptions of Theorem \ref{estimates},  Assume moreover that $\vb$ is controlled by $\vu$ : for every $t\in (0,T)$,
$$ \| \vb(t,.)\|_{L^3_{w_{3\gamma/2}}}\leq C_0 \|\vu(t,.)\|_{L^3_{w_{3\gamma/2}}}.$$
Then there exists a constant $C_\gamma\geq 1$ such that if $T_0<T$ is such that
$$ C_\gamma (1+C_0^4) \left(1+C_0^4+\|\vu_0\|_{L^2_{w_\gamma}}^2+\int_0^{T_0} \|\mathbb{F}\|_{L^2_{w_\gamma}}^2\, ds\right)^2\, T_0\leq 1$$ then
$$ \sup_{0\leq t\leq T_0} \|\ \vu(t,.)\|_{L^2_{w_\gamma}}^2 \leq
 C_\gamma (1+ C_0^4 + \|\vu_0\|_{L^2_{w_\gamma}}^2 +\int_0^{T_0} \|\mathbb{F}\|_{L^2_{w_\gamma}}^2\, ds
)$$ and 
$$  { \int_0^{T_0} \|\vN\vu\|_{L^2_{w_\gamma}}^2\, ds }\leq  
 C_\gamma (1+ C_0^4 + \|\vu_0\|_{L^2_{w_\gamma}}^2 +\int_0^{T_0} \|\mathbb{F}\|_{L^2_{w_\gamma}}^2\, ds).$$
\end{corollary}

\Proof We start from inequality (\ref{gronwb}) :  
 \begin{equation*}  \begin{split}  \|\vu(t,. )\|_{L^2_{w_\gamma}}^2 &+  \int_0^t  \|\vN\vu\|_{L^2_{w_\gamma}}^2\, ds  \\ \leq& \|\vu_0 \|_{L^2_{w_\gamma}}^2  +C_\gamma \int_0^t \|  \mathbb{F}(s,.)\|^2_{L^2_{w_\gamma}}\, ds + C_\gamma  \int_0^t  (1+  \|\vb(s,. )\|_{L^3_{w_{3\gamma/2}}}^2 )  \|\vu(s,. )\|_{L^2_{w_\gamma}}^2  \, ds\end{split}\end{equation*} 
 We write 
 $$  \|\vb(s,. )\|_{L^3_{w_{3\gamma/2}}}^2 \leq C_0^2  \|\vu(s,. )\|_{L^3_{w_{3\gamma/2}}}^2 \leq C_0^2 C_\gamma  \|u\|_{L^2_{w_\gamma}} (\|u\|_{L^2_{w_\gamma}}+\|\vN \vu\|_{L^2_{w_\gamma}}).$$
 This gives
  \begin{equation*}  \begin{split}  \|\vu(t,. ) & \|_{L^2_{w_\gamma}}^2  + \frac{1}{2}\int \| \vN \vu \|^2_{L^2_{w_\gamma}} \, ds  \\
   \leq &\|\vu_0 \|_{L^2_{w_\gamma}}^2  + C_\gamma \int_0^t \|  \mathbb{F}(s,.)\|^2_{L^2_{w_\gamma}}\, ds \\  &  + C_\gamma  \int_0^t   \|\vu(s,. )\|_{L^2_{w_\gamma}}^2  +   C_0^2  \|\vu(s,. )\|_{L^2_{w_\gamma}}^4 + C_0^4  \|\vu(s,. )\|_{L^2_{w_\gamma}}^6  \, ds
  \\  \leq &\|\vu_0 \|_{L^2_{w_\gamma}}^2  +C_\gamma \int_0^t \|  \mathbb{F}(s,.)\|^2_{L^2_{w_\gamma}}\, ds  +2 C_\gamma  \int_0^t   \|\vu(s,. )\|_{L^2_{w_\gamma}}^2  +    C_0^4  \|\vu(s,. )\|_{L^2_{w_\gamma}}^6  \, ds.  \end{split}\end{equation*}  For $t\leq T_0$, we get
  
  \begin{align*}
      \|\vu(t,. ) &  \|_{L^2_{w_\gamma}}^2  + \frac{1}{2}\int \| \vN \vu \|^2_{L^2_{w_\gamma}} \, ds  \\
      &\leq \|\vu_0 \|_{L^2_{w_\gamma}}^2  +C_\gamma \int_0^{T_0} \|\mathbb{F}\|_{L^2_{w_\gamma}}^2\, ds + C_\gamma (1+C_0^4) \int_0^t \|\vu(t,. )\|_{L^2_{w_\gamma}}^2+ \|\vu(t,. )\|_{L^2_{w_\gamma}}^6\, ds
  \end{align*}
  and we may conclude with Lemma \ref{gronwallnl}.\Endproof

\section{Stability of solutions for the advection-diffusion problem.}
  \subsection{The Rellich lemma.}
  
  We recall the Rellich lemma :
  \begin{lemma}[Rellich] If  $s>0$ and  $ (f_n)$ is a sequence of functions on $\mathbb{R}^d$ such that 
  \begin{itemize}
  \item[$\bullet$] the family $(f_n)$ is bounded in $H^s(\mathbb{R}^d)$
  \item[$\bullet$] there is a compact subset of $\mathbb{R}^d$ such that the support of each $f_n$ is included in $K$
\end{itemize}
then there exists a subsequence $(f_{n_k})$ such that $f_{n_k}$ is strongly convergent in $L^2(\mathbb{R}^d)$. 
\end{lemma}  
  
  We shall use a variant of this lemma (see \cite{LR02}) : 
  
  \begin{lemma}[space-time Rellich] \label{spacetime} If  $s>0$, $\sigma\in\mathbb{R}$ and  $ (f_n)$ is a sequence of functions on $(0,T)\times \mathbb{R}^d$ such that, for all $T_0\in (0,T)$ and all $\varphi\in\mathcal{D}(\mathbb{R}^3)$
  \begin{itemize}
  \item[$\bullet$]  $\varphi f_n$ is bounded in $L^2((0,T_0), H^s)$ 
  \item[$\bullet$]   $\varphi \partial_t f_n$ is bounded in $L^2((0,T_0), H^\sigma)$ 
\end{itemize}
then there exists a subsequence $(f_{n_k})$ such that $f_{n_k}$ is strongly convergent in $L^2_{\rm loc}([0,T)\times \mathbb{R}^3)$ : if $f_\infty$ is the limit, we have for all $T_0\in (0,T)$ and all $R_0>0$
$$\lim_{n_k\rightarrow +\infty} \int_0^{T_0} \int_{\vert x\vert\leq R} \vert f_{n_k}-f_\infty\vert^2 \, dx\, dt=0.$$
\end{lemma}  
\Proof    With no loss of generality, we may assume that $\sigma<\min(1,s)$.
 Define $g$ by $g_n(t,x)=\alpha(t)\varphi(x) f_n(t,x)$ if $t>0$ and  $g_n(t,x)=\alpha(t)\varphi(x) f_n(-t,x)$ if $t<0$, where $\alpha\in\mathcal{C}^\infty$ on $(0,T)$, is equal to $1$ on $[0,T_0]$ and equal to $0$ for $t> \frac{T+T_0}2$,  and $\varphi(x)=1$ on $B(0,R_0)$.
   Then the support of $g_n$ is contained in   $[-\frac{T+T_0}2,\frac{T+T_0}2]\times {\rm Supp} \,\varphi$.
  Moreover, $g_n$ is bounded in $L^2_t H^s$ and $\partial_t g_n$ is bounded in $L^2 H^\sigma$ so that $g_n $ is bounded in $H^\rho(\mathbb{R}\times\mathbb{R}^3)$ with $\rho= \frac s{s+1-\sigma}$ (just write $(1+\tau^2+\xi^2)^{\frac s{s+1-\sigma}} \leq \left((1+\tau^2)(1+\xi^2)^\sigma\right)^{\frac s{s+1-\sigma}} \left((1+\xi^2)^s\right)^{\frac{1-\sigma}{s+1-\sigma}}$).. By the Rellich lemma, we know that there is a subsequence $g_{n_k}$ which is strongly convergent in $L^2(\mathbb{R}\times\mathbb{R}^3)$, thus a subsequence $f_{n_k}$ which is strongly convergent in $L^2((0,T_0)\times B(0,R_0))$.
  
  We then iterate this argument for an increasing sequence of times $T_0<T_1<\dots <T_N\rightarrow T$ and an increasing sequence of radii $R_0<R_1<\dots <R_N\rightarrow +\infty$ and finish the proof. by the classical  diagonal process of Cantor.\Endproof
    \subsection{Proof of Theorem \ref{stability}.}

Assume that   $\vu_{0,n}$ is strongly convergent to $\vu_{0,\infty}$ in $L^2_{w_\gamma}$ and that the sequence $\mathbb{F}_n$ is strongly convergent to $\mathbb{F}_\infty$ in $L^2((0,T), L^2_{w_\gamma})$, and assume that   the sequence $\vb_n$ is bounded in $L^3((0,T),   L^3_{w_{3\gamma/2}})$. Then, by Theorem \ref{estimates} and Corollary \ref{passive},  we know that $\vu_n$ is bounded in $L^\infty((0,T), L^2_{w_\gamma})$ and $\vN\vu_n$ is bounded  in $L^2((0,T), L^2_{w_\gamma})$. In particular, writing
$p_n=p_{n,1}+ p_{n,2}$ with
 $$ p_{n,1}=\sum_{i=1}^3\sum_{j=1}^3 R_iR_j(b_{n,i}u_{n,j})  \text{ and }  p_{n,2}=-\sum_{i=1}^3\sum_{j=1}^3 R_iR_j( F_{n,i,j})$$
we get that  $p_{n,1}$ is bounded in   $L^{3}((0,T),L^{6/5}_{w_{\frac {6\gamma}5}})$ and $p_{n,2}$ is bounded in     $L^{2}((0,T),L^{2}_{w_\gamma})$.

If $\varphi\in \mathcal{D}(\mathbb{R}^3)$, we find that $\varphi \vu_n$ is bounded in $L^2((0,T), H^1)$
and, writing
$$ \partial_t \vu_n= \Delta \vu_n  - \left( \sum_{i=1}^3 \partial_i(b_{n,i}\vu_n) +\vN p_{n,1}\right)+ \left(\vN\cdot \mathbb{F}_n-\vN p_{n,2}\right), $$
$\varphi\partial_t\vu_n$ is bounded in $L^2 L^2 + L^2 W^{-1,6/5}+ L^2 H^{-1}\subset L^2((0,T), H^{-2})$. Thus, by Lemma \ref{spacetime},  there exists  $\vu_\infty$  and an increasing  sequence $(n_k)_{k\in\mathbb{N}}$ with values in $\mathbb{N}$ such that
 $\vu_{n_k}$ converges strongly  to $\vu_\infty$ in  $L^2_{\rm loc}([0,T)\times\mathbb{R}^3)$ : for every  $T_0\in (0,T)$ and every $R>0$, we have
  $$\lim_{k\rightarrow +\infty} \int_0^{T_0} \int_{\vert y\vert<R} \vert \vu_{n_k}(s,y)-\vu_\infty(s,y)\vert^2\, dy\, ds=0.$$
As $\vu_n$ is bounded in $L^\infty((0,T), L^2_{w_\gamma})$ and $\vN\vu_n$ is bounded  in $L^2((0,T), L^2_{w_\gamma})$, the convergence of $\vu_{n_k}$ to $\vu_\infty$ in $\mathcal{D}'((0,T)\times \mathbb{R}^3)$ implies that $\vu_{n_k}$ converges *-weakly to $\vu_\infty$ in $L^\infty((0,T), L^2_{w_\gamma})$ and $\vN\vu_{n_k}$ converges weakly to $\vN\vu_\infty$ in $L^2((0,T),L^2_{w_\gamma})$.

By Banach--Alaoglu's theorem, we may assume that   there exists   $\vb_\infty$  such that $\vb_{n_k}$ converges weakly to $\vb_\infty$ in $L^3((0,T), L^3_{w_{3\gamma/2}})$.
In particular $b_{n_k,i} u_{n_k,j}$ is weakly convergent in $(L^{6/5}L^{6/5})_{\rm loc}$ and thus in $\mathcal{D}'((0,T)\times \mathbb{R}^3)$; as it is bounded in  $L^{3}((0,T),L^{6/5}_{w_{\frac {6\gamma}5}})$, it is weakly convergent in  $L^{3}((0,T),L^{6/5}_{w_{\frac {6\gamma}5}})$ to $b_{\infty,i}u_{\infty,j}$. Let  $$p_{\infty,1}=\sum_{i=1}^3\sum_{j=1}^3 R_i R_j(b_{\infty,i}u_{\infty,j}) \text{ and  } p_{\infty,2}=-\sum_{i=1}^3\sum_{j=1}^3 R_i R_j( F_{\infty,i,j}).$$ As the Riesz transforms are bounded on $L^{6/5}_{w_{\frac {6\gamma}5}}$ and on $L^2_{w_\gamma}$, we find that $p_{n_k,1}$ is weakly convergent in  $L^{3}((0,T),L^{6/5}_{w_{\frac {6\gamma}5}})$  to $p_{\infty,1}$ and that $p_{n_k,2}$  is strongly convergent in  $L^2((0,T),L^2_{w_\gamma})$  to $p_{\infty,2}$.

In particular, we find that in $\mathcal{D}'((0,T)\times\mathbb{R}^3)$
$$ \partial_t \vu_\infty=\Delta \vu_\infty-\sum_{i=1}^3 \partial_i (b_{\infty,i}\vu_\infty)-\vN(p_{\infty,1}+p_{\infty,2})+\vN\cdot \mathbb{F}_\infty.$$
In particular, $\partial_t\vu_\infty$ is locally in $L^2 H^{-2}$, and thus $\vu_\infty$ has representative such that  $t\mapsto \vu_\infty(t,.)$ is continuous from $[0,T)$ to $\mathcal{D}'(\mathbb{R}^3)$ and coincides with $\vu_\infty(0,.)+\int_0^t \partial_t \vu_\infty\, ds$. In $\mathcal{D}'((0,T)\times \mathbb{R}^3)$, we have that $$\vu_\infty(0,.)+\int_0^t \partial_t \vu_\infty\, ds=\vu_\infty=\lim_{n_k\rightarrow +\infty} \vu_{n_k}=\lim_{n_k\rightarrow +\infty}
\vu_{0,n_k}+ \int_0^t \partial_t \vu_{n_k}\, ds=\vu_{0, \infty}+\int_0^t \partial_t\vu_\infty\, ds$$
Thus, $\vu_\infty(0,.)=\vu_{0,\infty}$, and $\vu_\infty$ is a solution of $(AD_\infty)$.

 Next, we define 
 $$ A_n=  - \partial_t(\frac {\vert\vu_n\vert^2}2)+\Delta(\frac {\vert\vu_n\vert^2}2)-\vN\cdot\left( \frac{\vert\vu_n\vert^2}2\vb_n\right)-\vN\cdot(p_n\vu_n) + \vu_n\cdot(\vN\cdot\mathbb{F}_n) =\vert\vN \vu_n\vert^2 +\mu_n.$$
As $\vu_n$ is bounded in $L^\infty((0,T), L^2_{w_\gamma})$ and $\vN\vu_n$ is bounded  in $L^2((0,T), L^2_{w_\gamma})$, it is bounded in $L^2((0,T), L^6_{w_{3\gamma/2}})$ and by interpolation with $L^\infty((0,T), L^2_{w_\gamma})$  it is bounded in $L^{10/3}((0,T), L^{10/3}_{w_{5\gamma/3}})$. Thus, $u_{n_k}$ is locally bounded in $L^{10/3}L^{10/3}$ and locally strongly convergent in $L^2 L^2$; it is then strongly convergent in $L^3 L^3$.  Thus, $A_{n_k}$ is convergent in $\mathcal{D}'((0,T)\times\mathbb{R}^3)$ to
 $$ A_\infty=  - \partial_t(\frac {\vert\vu_\infty\vert^2}2)+\Delta(\frac {\vert\vu_\infty\vert^2}2)-\vN\cdot\left( \frac{\vert\vu_\infty\vert^2}2\vb_\infty\right)-\vN\cdot(p_\infty\vu_\infty) + \vu_\infty\cdot(\vN\cdot\mathbb{F}_\infty) .$$ 
 In particular, $A_\infty =\lim_{n_k\rightarrow +\infty} \vert\vN \vu_{n_k}\vert^2 +\mu_{n_k}$. If $\Phi\in \mathcal{D}((0,T)\times \mathbb{R}^3)$ is non-negative, we have
$$ \iint A_\infty \Phi\, dx\, ds=\lim_{n_k\rightarrow +\infty}\iint A_{n_k} \Phi\, dx\, ds\geq \limsup_{n_k\rightarrow +\infty}\iint   \vert\vN \vu_{n_k}\vert^2 \Phi\, dx\, ds \geq  \iint   \vert\vN \vu_{\infty}\vert^2 \Phi\, dx\, ds$$ (since $\sqrt\Phi \vN \vu_{n_k}$ is weakly convergent to $\sqrt \Phi \vN\vu_\infty$ in $L^2L^2$). Thus, there   exists a non-negative locally finite measure $\mu_\infty$ on $(0,T)\times\mathbb{R}^3$ such that $A_\infty=\vert\vN \vu_\infty\vert^2 +\mu_\infty$, i.e. such that
 $$ \partial_t(\frac {\vert\vu_\infty\vert^2}2)=\Delta(\frac {\vert\vu_\infty\vert^2}2)-\vert\vN \vu_\infty\vert^2- \vN\cdot\left( \frac{\vert\vu_\infty\vert^2}2\vb_\infty\right)-\vN\cdot(p_\infty\vu_\infty) + \vu\cdot(\vN\cdot\mathbb{F}_\infty)-\mu_\infty.$$

Finally, we start from inequality  (\ref{localineq}) :

 \begin{equation*}  \begin{split}
\int \frac{\vert \vu_n(t,x)\vert^2}2 \phi_R w_{\gamma,\epsilon}\, dx    & 
 \leq \int \frac{\vert \vu_{0,n}(x)\vert^2}2 \phi_R w_{\gamma,\epsilon}\, dx     
 \\ &-\sum_{i=1}^3 \int_0^t \int   \partial_i\vu_n\cdot \vu_n\,  (w_{\gamma,\epsilon}\partial_i \phi_R+\phi_R \partial_iw_{\gamma,\epsilon})\, dx\, ds
\\&- \int_0^t  \int  \vert\vN\vu_n\vert^2\, \,  \phi_R w_{\gamma,\epsilon} dx\, ds
\\&  + \sum_{i=1}^3 \int_0^t \int  \frac{\vert\vu_n\vert^2}2 b_{n,i}     (w_{\gamma,\epsilon}\partial_i \phi_R+\phi_R \partial_iw_{\gamma,\epsilon})\, dx\, ds
\\& +\sum_{i=1}^3 \int_0^t    \int     p_nu_{n,i} (w_{\gamma,\epsilon}\partial_i \phi_R+\phi_R \partial_iw_{\gamma,\epsilon})\, dx\, ds
\\&- \sum_{i=1}^3\sum_{j=1}^3 \int_0^t \int    F_{n,i,j} u_{n,j}    (w_{\gamma,\epsilon}\partial_i \phi_R+\phi_R \partial_iw_{\gamma,\epsilon})\, dx\, ds
\\&  - \sum_{i=1}^3 \sum_{j=1}^3\int_0^t  \int    F_{n,i,j}\partial_iu_{n,}\    \phi_R w_{\gamma,\epsilon}\, dx\, ds.
\end{split}\end{equation*}
This gives  
 \begin{equation*}  \begin{split}
\limsup_{n_k\rightarrow +\infty}\int \frac{\vert \vu_{n_k}(t,x)\vert^2}2 \phi_R w_{\gamma,\epsilon}\, dx&     
+\int_0^t  \int  \vert\vN\vu_{n_k}\vert^2\, \,  \phi_R w_{\gamma,\epsilon} dx\, ds \\ \leq& \int \frac{\vert \vu_{0,\infty}(x)\vert^2}2 \phi_R w_{\gamma,\epsilon}\, dx     
 \\ &-\sum_{i=1}^3 \int_0^t \int   \partial_i\vu_\infty\cdot \vu_\infty\,  (w_{\gamma,\epsilon}\partial_i \phi_R+\phi_R \partial_iw_{\gamma,\epsilon})\, dx\, ds
\\&  + \sum_{i=1}^3 \int_0^t \int  \frac{\vert\vu_\infty\vert^2}2 b_{\infty,i}     (w_{\gamma,\epsilon}\partial_i \phi_R+\phi_R \partial_iw_{\gamma,\epsilon})\, dx\, ds
\\& +\sum_{i=1}^3 \int_0^t    \int     p_\infty u_{\infty,i} (w_{\gamma,\epsilon}\partial_i \phi_R+\phi_R \partial_iw_{\gamma,\epsilon})\, dx\, ds
\\&- \sum_{i=1}^3\sum_{j=1}^3 \int_0^t \int    F_{\infty,i,j} u_{\infty,j}    (w_{\gamma,\epsilon}\partial_i \phi_R+\phi_R \partial_iw_{\gamma,\epsilon})\, dx\, ds
\\&  - \sum_{i=1}^3 \sum_{j=1}^3\int_0^t  \int    F_{\infty,i,j}\partial_iu_{\infty,j}\    \phi_R w_{\gamma,\epsilon}\, dx\, ds.
\end{split}\end{equation*}

As we have $$ \vu_{n_k}= 
\vu_{0, n_k}+ \int_0^t \partial_t \vu_{n_k}\, ds $$ we see that $\vu_{n_k}(t,.)$ is convergent to $\vu_\infty(t,.)$ in $\mathcal{D}'(\mathbb{R}^3)$, hence is weakly convergent in $L^2_{\rm loc}$ (as it is bounded in $L^2_{w_\gamma}$), so that :
$$
 \int \frac{\vert \vu_\infty(t,x)\vert^2}2 \phi_R w_{\gamma,\epsilon}\, dx   \leq  
\limsup_{n_k\rightarrow +\infty}\int \frac{\vert \vu_{n_k}(t,x)\vert^2}2 \phi_R w_{\gamma,\epsilon}\, dx .$$    
  Similarly, as $\vN\vu_{n_k}$ is weakly convergent in $L^2L^2_{w_\gamma}$, we have
  $$
\int_0^t \int \frac{\vert \vN\vu_\infty(s,x)\vert^2}2 \phi_R w_{\gamma,\epsilon}\, dx \, ds  \leq  
\limsup_{n_k\rightarrow +\infty} \int_0^t\int \frac{\vert \vN\vu_{n_k}(s,x)\vert^2}2 \phi_R w_{\gamma,\epsilon}\, dx \, ds.$$    

Thus, letting $R$ go to $+\infty$ and then $\epsilon$ go to $0$, we find by dominated convergence that,
  for every $t\in (0,T)$, we have
 \begin{equation*}\begin{split}
 \|\vu_\infty(t,. )\|_{L^2_{w_\gamma}}^2 &+2\int_0^t \|\vN \vu_\infty(s,.)\|_{L^2_{w_\gamma}}^2\, ds
 \\  \leq& \|\vu_{0, \infty} \|_{L^2_{w_\gamma}}^2 -\int_0^t \int   \vN\vert\vu_\infty\vert^2\cdot\vN w_\gamma\, dx\, ds+\int_0^t\int (\vert\vu_\infty\vert^2 \vb_\infty+2p_\infty\vu_\infty) \cdot \vN(w_\gamma)\, dx\, ds\\ &\phantom{space space} -2\sum_{i=1}^3\sum_{j=1}^3 \int_0^t\int F_{\infty,i,j} (\partial_i u_{\infty,j}) w_\gamma +F_{\infty,i,j} u_{\infty,i} \partial_j(w_\gamma)\, dx\, ds.
 \end{split}\end{equation*} Letting $t$ go to $0$, we find
 $$\limsup_{t\rightarrow 0}  \|\vu_\infty(t,. )\|_{L^2_{w_\gamma}}^2\leq  \|\vu_{0,\infty} \|_{L^2_{w_\gamma}}^2 .$$   On the other hand, we know that $\vu_\infty$ is weakly continuous in $L^2_{w_\gamma}$ and thus we have 
  $$    \|\vu_{0,\infty} \|_{L^2_{w_\gamma}}^2 \leq   \liminf_{t\rightarrow 0}  \|\vu_\infty(t,. )\|_{L^2_{w_\gamma}}^2 .$$   
  This gives $    \|\vu_{0,\infty} \|_{L^2_{w_\gamma}}^2 =   \lim_{t\rightarrow 0}  \|\vu_\infty(t,. )\|_{L^2_{w_\gamma}}^2 $, which allows to turn the weak convergence into a strong convergence. Theorem \ref{stability} is proven. \Endproof

\section{Solutions  of the Navier--Stokes problem with initial data in  $L^2_{w_\gamma}$.}

We now prove Theorem \ref{weightedNS}.  The idea is to approximate the problem by a Navier--Stokes problem in $L^2$, then use the a priori estimates (Theorem \ref{estimates}) and the stability theorem (Theorem \ref{stability}) to find a solution to the Navier--Stokes problem with data in $L^2_{w_\gamma})$.

\subsection{Approximation by square integrable data.}

\begin{lemma}[Leray's projection operator]\label{projec}
Let $0<\delta<3$ and $1<r<+\infty$. If $\vv$ is a vector field on $\mathbb{R}^3$ such that $\vv\in L^r_{w_{\delta}}$, then there exists a unique decompostion
$$ \vv=\vv_\sigma+\vv_\nabla$$ 
such that
\begin{itemize}
\item[$\bullet$] $\vv_\sigma\in L^r_{w_\delta}$ and $\vN\cdot \vv_\sigma=0$.
\item[$\bullet$] $\vv_\nabla\in L^r_{w_\delta}$ and $\vN\wedge \vv_\nabla=0$.
\end{itemize}
We shall write $\vv_\sigma=\mathbb{P}\vv$, where $\mathbb{P}$ is Leray's projection operator.

Similarly, if $\vv$ is a distribution vector field of the type $\vv=\vN\cdot\mathbb{G}$ with $\mathbb{G}\in L^r_{w_\delta}$ then  there exists a unique decompostion
$$ \vv=\vv_\sigma+\vv_\nabla$$ 
such that
\begin{itemize}
\item[$\bullet$]  there exists $\mathbb{H}\in L^r_{w_\delta}$ such that $\vv_\sigma\ = \vN\cdot \mathbb{H}$  and $\vN\cdot \vv_\sigma=0$.
\item[$\bullet$] there exists $q\in L^r_{w_\delta}$ such that $\vv_\nabla =\vN q$ (and thus  $\vN\wedge \vv_\nabla=0$).
\end{itemize}
We shall still write  $\vv_\sigma=\mathbb{P}\vv$.  Moreover, the function $q$ is given by
$$ q=- \sum_{i=1}^3\sum_{j=1}^3 R_i R_j (G_{i,j}). $$
\end{lemma}

\Proof As $w_\delta\in \mathcal{A}_r$ the Riesz transforms are bounded on $L^r_{w_\delta}$. 
Using the identity
$$ \Delta \vv=\vN(\vN\cdot \vv)-\vN\wedge(\vN\wedge\vv)$$ we find (if the decomposition exists) that  
$$ \Delta \vv_\sigma = -\vN\wedge(\vN\wedge\vv_\sigma) =-\vN\wedge(\vN\wedge\vv) \text{ and } \Delta\vv_\nabla = \vN(\vN\cdot \vv_\nabla)=\vN(\vN\cdot \vv).$$

This proves the uniqueness. By linearity, we just have to prove that $\vv=0\implies \vv_\nabla=0$. We have $\Delta\vv_\nabla=0$, and thus $\vv_\nabla$ is harmonic; as it belongs to $\mathcal{S}'$, we  find that it is a polynomial. But a polynomial which belongs to $L^r_{w_\delta}$ must be equal to $0$. 
Similarly, if $\vv_\nabla=\vN q$, then $\Delta q=\vN\cdot \vv_\nabla=\vN\cdot \vv=0$; thus $q$ is harmonic and belongs to $L^r_{w_\delta}$, hence $q=0$.

For the existence, it is enough to check that  $v_{\nabla,i} =-\sum_{j=1}^3 R_iR_j v_j $ in the first case and $\vv_\nabla=\vN q$ with $q=  \sum_{i=1}^3\sum_{j=1}^3 R_i R_j (G_{i,j})$ in the second case fulfill the conclusions of the lemma.\Endproof

\begin{lemma}\label{approx}  Let $0<\gamma\leq 2$. Let $\vu_{0}$ be a divergence-free vector field such that $\vu_0\in L^2_{w_\gamma}(\mathbb{R}^3)$ and $\mathbb{F}$ be a tensor $\mathbb{F}(t,x)=\left(F_{i,j}(t,x)\right)_{1\leq i,j\leq 3}$ such that $\mathbb{F}\in L^2((0,+\infty), L^2_{w_\gamma})$.  Let  $\phi\in\mathcal{D}(\mathbb{R}^3)$  be  a non-negative  function which is equal to $1$ for $\vert x\vert\leq 1$ and to $0$ for $\vert x\vert\geq 2$. For $R>0$, we define $\phi_R(x)=\phi(\frac x R)$, $\vu_{0,R}=\mathbb{P}(\phi_R \vu_0)$ and $\mathbb{F}_R=\phi_R \mathbb{F}$.  Then $\vu_{0,R}$ is a divergence-free square integrable vector field and $\lim_{R\rightarrow +\infty} \|\vu_{0,R}-\vu_0\|_{L^2_{w_\gamma}}=0$. Similarly, $\mathbb{F}_R$ belongs to $L^2 L^2$ and  $\lim_{R\rightarrow +\infty} \|\mathbb{F}_R-\mathbb{F}\|_{L^2((0,+\infty),L^2_{w_\gamma})}=0$. 
\end{lemma}

 \Proof  By dominated convergence, we have $\lim_{R\rightarrow +\infty} \|\phi_R \vu_0-\vu_0\|_{L^2_{w_\gamma}}=0$. We conclude by writing $\vu_{0,R}-\vu_0=\mathbb{P}(\phi_R\vu_0-\vu_0)$.
 \Endproof
 \subsection{Leray's mollification.}
 We want to solve the Navier--Stokes equations with initial value $\vu_0$ :
 \begin{equation*}  (NS) \left\{ \begin{matrix} \partial_t \vu= \Delta \vu  -(\vu\cdot \vN)\vu- \vN p +\vN\cdot \mathbb{F} \cr \cr \vN\cdot \vu=0,    \phantom{space space} \vu(0,.)=\vu_0
 \end{matrix}\right.\end{equation*}
 We begin with Leray's method \cite{Le34} for solving the problem in $L^2$ :  \begin{equation*}  (NS_R) \left\{ \begin{matrix} \partial_t \vu_R= \Delta \vu_R  -(\vu_R\cdot \vN)\vu_R- \vN p_R +\vN\cdot \mathbb{F_R} \cr \cr \vN\cdot \vu_R=0,    \phantom{space space} \vu_R(0,.)=\vu_{0,R}
 \end{matrix}\right.\end{equation*}
The idea of Leray is to mollify the non-linearity by replacing $\vu_R\cdot\vN$ by $(\vu_R*\theta_\epsilon)\cdot \vN$, where $\theta(x)=\frac 1{\epsilon^3}\theta(\frac x\epsilon)$, $\theta\in\mathcal{D}(\mathbb{R}^3)$, $\theta$ is non-negative and radially decreasing and $\int\theta\, dx=1$. We thus solve the problem 
 \begin{equation*}  (NS_{R,\epsilon}) \left\{ \begin{matrix} \partial_t \vu_{R,\epsilon}= \Delta \vu_{R,\epsilon}  -((\vu_{R,\epsilon}*\theta_\epsilon)\cdot \vN)\vu_{R,\epsilon}- \vN p_{R,\epsilon} +\vN\cdot \mathbb{F}_{R} \cr \cr \vN\cdot \vu_{R,\epsilon}=0,    \phantom{space space} \vu_{R,\epsilon}(0,.)=\vu_{0,R}
 \end{matrix}\right.\end{equation*}

The classical result of Leray states that the problem $(NS_{R,\epsilon})$ is well-posed~:

\begin{lemma}\label{leray} Let $\vv_0\in L^2$ be a divergence-free vector field. Let $\mathbb{G}\in L^2((0,+\infty), L^2)$. Then  
 the problem 
 \begin{equation*}  (NS_{\epsilon}) \left\{ \begin{matrix} \partial_t \vv_\epsilon = \Delta \vv_\epsilon   -((\vv_\epsilon *\theta_\epsilon)\cdot \vN)\vv_\epsilon - \vN q_\epsilon  +\vN\cdot \mathbb{G}  \cr \cr \vN\cdot \vv_\epsilon =0,    \phantom{space space} \vv_\epsilon (0,.)=\vv_0
 \end{matrix}\right.\end{equation*}
 has a unique solution $\vv_\epsilon$ in $L^\infty((0,+\infty), L^2)\cap L^2((0,+\infty),\dot H^1)$. Moreover, this solution belongs to $\mathcal{C}([0,+\infty), L^2)$.
\end{lemma}

 \subsection{Proof of Theorem \ref{weightedNS} (local existence)}
 
 We use Lemma \ref{leray}  and find a solution $\vu_{R,\epsilon}$ to the problem $(NS_{R,\epsilon})$. 
Then we check that $\vu_{R,\epsilon}$ fulfills the assumptions of Theorem  \ref{estimates} and of Corollary \ref{active} :
\begin{itemize}
 \item[$\bullet$]  $\vu_{R,\epsilon}$ belongs to $L^\infty((0,T), L^2_{w_\gamma})$ and $\vN\vu_{R,\epsilon}$ belongs to $L^2((0,T),L^2_{w_\gamma})$  \item[$\bullet$] the map $t\in [0,+\infty)\mapsto \vu_{R, \epsilon} (t,.)$ is weakly continuous from $[0,+\infty)$ to $L^2_{w_\gamma}$, and is strongly continuous at $t=0$ :
 $$ \lim_{t\rightarrow 0} \|\vu_{R,\epsilon} (t,.)-\vu_{0, R}\|_{L^2_{w_\gamma}}=0.$$
 \item[$\bullet$]   on $(0,T)\times\mathbb{R}^3$, $\vu_{R,\epsilon}$ fulfills the energy equality : 
 \begin{equation*} \partial_t(\frac {\vert\vu_{R,\epsilon}\vert^2}2)=\Delta(\frac {\vert\vu_{R,\epsilon}\vert^2}2)-\vert\vN \vu_{R,\epsilon}\vert^2- \vN\cdot\left( \frac{\vert\vu\vert^2}2\vb_{R,\epsilon}\right)-\vN\cdot(p_{R,\epsilon}\vu_{R,\epsilon}) + \vu_{R,\epsilon}\cdot(\vN\cdot\mathbb{F}_R).\end{equation*}
with $\vb_{R,\epsilon}=\vu_{R,\epsilon}*\theta_\epsilon$.\item[$\bullet$] $\vb_{R,\epsilon}$ is controlled by $\vu_{R,\epsilon}$ : for every $t\in (0,T)$,
$$ \| \vb_{R,\epsilon}(t,.)\|_{L^3_{w_{3\gamma/2}}}\leq \|\mathcal{M}_{\vu_{R,\epsilon}(t,.)}\|_{L^3_{w_{3\gamma/2}}} \leq C_0 \|\vu_{R,\epsilon}(t,.)\|_{L^3_{w_{3\gamma/2}}}.$$
\end{itemize}
Thus, we know that, for every time $T_0$ such that
$$ C_\gamma (1+C_0^4) \left(1+C_0^4+\|\vu_{0,R}\|_{L^2_{w_\gamma}}^2+\int_0^{T_0} \|\mathbb{F}_R\|_{L^2_{w_\gamma}}^2\, ds\right)^2\, T_0\leq 1$$ we have
$$ \sup_{0\leq t\leq T_0} \|\ \vu_{R,\epsilon}(t,.)\|_{L^2_{w_\gamma}}^2 \leq
 C_\gamma (1+ C_0^4 + \|\vu_{0,R}\|_{L^2_{w_\gamma}}^2 +\int_0^{T_0} \|\mathbb{F}_R\|_{L^2_{w_\gamma}}^2\, ds
)$$ and 
$$  { \int_0^{T_0} \|\vN\vu_{R,\epsilon}\|_{L^2_{w_\gamma}}^2\, ds }\leq  
 C_\gamma (1+ C_0^4 + \|\vu_{0,R}\|_{L^2_{w_\gamma}}^2 +\int_0^{T_0} \|\mathbb{F}_R\|_{L^2_{w_\gamma}}^2\, ds).$$
Moreover, we have that 
$$  \|\vu_{0,R}\|_{L^2_{w_\gamma}}\leq C_\gamma  \|\vu_0\|_{L^2_{w_\gamma}}  \text{ and }  \|\mathbb{F}_R\|_{L^2_{w_\gamma}} \leq \|\mathbb{F}\|_{L^2_{w_\gamma}}$$
  so that
\begin{equation*} \begin{split}  \|\vb_{R,\epsilon}\|_{L^3((0,T_0), L^3_{w_{3\gamma/2}}}\leq &C_\gamma   \|\vu_{R,\epsilon}\|_{L^3((0,T_0), L^3_{w_{3\gamma/2}}} \\ \leq& C'_\gamma T_0^{ \frac 1{12}} \left ( (1+\sqrt{ T_0}) \|\vu_{R,\epsilon}\|_{L^\infty((0,T_0), L^2_{w_\gamma})} + \|\vN\vu_{R,\epsilon}\|_{L^2((0,T_0), L^2_{w_\gamma})}\right)
\\  \leq&
 C''_\gamma  \sqrt{1+ C_0^4 + \|\vu_{0}\|_{L^2_{w_\gamma}}^2 +\int_0^{T_0} \|\mathbb{F}\|_{L^2_{w_\gamma}}^2\, ds}.
 \end{split}\end{equation*}
 
 Let $R_n\rightarrow +\infty$ and $\epsilon_n\rightarrow 0$. Let $\vu_{0,n}=\vu_{0,R_n}$, $\mathbb{F}_n=\mathbb{F}_{R_n}$, $\vb_n=\vb_{R_n,\epsilon_n}$ and $\vu_n=\vu_{R_n,\epsilon_n}$.  We may then apply Theorem \ref{stability}, since  $\vu_{0,n}$ is strongly convergent to $\vu_{0}$ in $L^2_{w_\gamma}$,   $\mathbb{F}_n$ is strongly convergent to $\mathbb{F}$  in $L^2((0,T_0), L^2_{w_\gamma})$,  and   the sequence $\vb_n$ is bounded in $L^3((0,T_0),   L^3_{w_{3\gamma/2}})$. Thus  there exists $p$, $\vu$,    $\vb$   and an increasing  sequence $(n_k)_{k\in\mathbb{N}}$ with values in $\mathbb{N}$ such that

  \begin{itemize} 
 \item[$\bullet$] $\vu_{n_k}$ converges *-weakly to $\vu$ in $L^\infty((0,T_0), L^2_{w_\gamma})$, $\vN\vu_{n_k}$ converges weakly to $\vN\vu$ in $L^2((0,T_0),L^2_{w_\gamma})$
 \item[$\bullet$] $\vb_{n_k}$ converges weakly to $\vb$ in $L^3((0,T_0), L^3_{w_{3\gamma/2}})$, $p_{n_k}$ converges weakly to $p$ in $L^{3}((0,T_0),L^{6/5}_{w_{\frac {6\gamma}5}})+L^{2}((0,T_0),L^{2}_{w_\gamma})$
  \item[$\bullet$] $\vu_{n_k}$ converges strongly  to $\vu$ in  $L^2_{\rm loc}([0,T_0)\times\mathbb{R}^3)$.  \end{itemize}

 Moreover, $\vu$ is  a solution of the   advection-diffusion problem 
 \begin{equation*}   \left\{ \begin{matrix} \partial_t \vu= \Delta \vu  -(\vb\cdot \vN)\vu- \vN p +\vN\cdot \mathbb{F} \cr \cr \vN\cdot \vu=0,    \phantom{space space} \vu(0,.)=\vu_{0}
 \end{matrix}\right.\end{equation*}
 and is
  such that :
 \begin{itemize}  
 \item[$\bullet$] the map $t\in [0,T_0)\mapsto \vu(t,.)$ is weakly continuous from $[0,T_0)$ to $L^2_{w_\gamma}$, and is strongly continuous at $t=0$ :
 $$ \lim_{t\rightarrow 0} \|\vu(t,.)-\vu_{0}\|_{L^2_{w_\gamma}}=0.$$
 \item[$\bullet$]  there exists a non-negative locally finite measure $\mu$ on $(0,T_0)\times\mathbb{R}^3$ such that
 $$ \partial_t(\frac {\vert\vu\vert^2}2)=\Delta(\frac {\vert\vu\vert^2}2)-\vert\vN \vu\vert^2- \vN\cdot\left( \frac{\vert\vu\vert^2}2\vb\right)-\vN\cdot(p \vu) + \vu\cdot(\vN\cdot\mathbb{F})-\mu.$$
 \end{itemize}

 Finally, as $\vb_n=\theta_{\epsilon_n}*(\vu_n-\vu)+ \theta_{\epsilon_n}*\vu$, we see that $\vb_{n_k}$ is strongly   convergent to $\vu$ in   $L^3_{\rm loc}([0,T_0)\times\mathbb{R}^3)$, so that $\vb=\vu$ : thus,  $\vu$ is a solution of the Navier--Stokes problem on $(0,T_0)$. (It is easy to check that   $$ p=\sum_{i=1}^3\sum_{j=1}^3 R_iR_j(u_iu_j-F_{i,j})$$ as $u_{i,n_k} u_{j,n_k}$ is weakly convergent to $u_i u_j$ in 
  $L^{4}((0,T_0),L^{6/5}_{w_{\frac {6\gamma}5}})$ and $w_{\frac {6\gamma}5}
\in \mathcal{A}_{6/5}$).

  \subsection{Proof of Theorem \ref{weightedNS} (global existence)}

 In order to finish the proof, we shall use the scaling properties of the Navier--Stokes equations : if $\lambda>0$, then $\vu$ is a solution of the Cauchy initial value problem for the Navier--Stokes equations on $(0,T)$ with initial value $\vu_0$ and forcing tensor $\mathbb{F}$  if and only if $\vu_\lambda(t,x)=\lambda \vu(\lambda^2 t,\lambda x)$ is a solution of the Navier--Stokes equations on $(0,T/\lambda^2)$ with initial value $\vu_{0,\lambda}(x)=\lambda \vu_0(\lambda x)$ and forcing tensor $\mathbb{F}_\lambda(t,x)=\lambda^2 \mathbb{F}(\lambda^2  t,\lambda x)$.
 
 We  take $\lambda>1$ and for $n\in\mathbb{N}$ we consider the Navier--Stokes problem with initial value $\vv_{0,n}=\lambda^n \vu_0(\lambda^n \cdot)$ and forcing tensor $\mathbb{F}_n=\lambda^{2n}\mathbb{F}(\lambda^{2n} \cdot , \lambda^n \cdot )$. Then we have seen that we can find a solution $\vv_n$ on $(0,T_n)$, with 
 $$   C_\gamma \left(1+\|\vv_{0,n}\|_{L^2_{w_\gamma}}^2+\int_0^{+\infty} \|\mathbb{F}_n\|_{L^2_{w_\gamma}}^2\, ds\right)^2\, T_n =  1.$$  
 Of course, we have $\vv_n(t,x)=\lambda^{n}\vu_n(\lambda^{2n}t,\lambda^n x)$ where $\vu_n$ is a solution of the   Navier--Stokes equations on $(0,\lambda^{2n}T_n)$ with initial value $\vu_0$ and forcing tensor $\mathbb{F}$
 
 \begin{lemma}
 $$\lim_{n\rightarrow +\infty} \frac{\lambda^n}{ 1+\|\vv_{0,n}\|_{L^2_{w_\gamma}}^2+\int_0^{+\infty} \|\mathbb{F}_n\|_{L^2_{w_\gamma}}^2\, ds}=+\infty.
 $$\end{lemma}
 
 \Proof We have
 $$ \|\vv_{0,n}\|_{L^2_{w_\gamma}}^2=\int \vert\vu_0(x)\vert^2 \lambda^{n(\gamma-1)} \frac{(1+\vert x\vert)^\gamma}{(\lambda^n+\vert x\vert)^\gamma}  w_\gamma(x)\, dx. $$ We have $$\lambda^{n(\gamma-1)}\leq \lambda^n$$ as $\gamma\leq 2$ and we have, by dominated convergence,
 $$\lim_{n\rightarrow +\infty}  \int \vert\vu_0(x)\vert^2   \frac{(1+\vert x\vert)^\gamma}{(\lambda^n+\vert x\vert)^\gamma}  w_\gamma(x)\, dx=0. $$ Similarly, we have 
 $$\int_0^{+\infty} \| \mathbb{F}_{n}\|_{L^2_{w_\gamma}}^2\, ds =\int_0^{+\infty}\int \vert \mathbb{F}(s,x)\vert^2 \lambda^{n(\gamma-1)} \frac{(1+\vert x\vert)^\gamma}{(\lambda^n+\vert x\vert)^\gamma}  w_\gamma(x)\, dx\, ds=o(\lambda^n). $$  Thus, $\lim_{n\rightarrow +\infty} \lambda^{2n} T_n=+\infty$.
 \Endproof

 Now, for a given $T>0$, if $  \lambda^{2n}T_n>T$ for $n\geq n_T$, then $\vu_n$ is a solution of the Navier-Stokes problem on $(0,T)$. Let $\vw_n(t,x)=\lambda^{n_T} \vu_n(\lambda^{2n_T}t,  \lambda^{n_T}x)$.  For $n\geq n_T$, $\vw_n$ is a solution of the Navier-Stokes problem on $(0,\lambda^{-2n_T} T)$ with initial value $\vv_{0,n_T}$ and forcing tensor $\mathbb{F}_{n_T}$. 
 As $\lambda^{-2n_T}T\leq T_{n_T}$, we have
  $$   C_\gamma \left(1+\|\vv_{0,n_T}\|_{L^2_{w_\gamma}}^2+\int_0^{+\infty} \|\mathbb{F}_{n_T}\|_{L^2_{w_\gamma}}^2\, ds\right)^2\, \lambda^{-2n_T} T \leq  1.$$  
By corollary \ref{active}, we have
$$ \sup_{0\leq t\leq   \lambda^{-2n_T} T} \|\ \vw_n(t,.)\|_{L^2_{w_\gamma}}^2 \leq
 C_\gamma (1+  \|\vv_{0,n_T}\|_{L^2_{w_\gamma}}^2 +\int_0^{\lambda^{-2n_T}T} \|\mathbb{F}_{n_T}\|_{L^2_{w_\gamma}}^2\, ds 
)$$ and 
$$  { \int_0^{\lambda^{-2n_T}T} \|\vN\vw_n\|_{L^2_{w_\gamma}}^2\, ds }\leq  
 C_\gamma (1+   \|\vv_{0,n_T}\|_{L^2_{w_\gamma}}^2 +\int_0^{\lambda^{-2n_T}T} \|\mathbb{F}_{n_T}\|_{L^2_{w_\gamma}}^2\, ds).$$
 We have
 $$ \|\vw_n\|_{L^2_{w_\gamma}}^2=\int \vert\vu_n(\lambda^{2n_T}t,x)\vert^2 \lambda^{n_T(\gamma-1)} \frac{(1+\vert x\vert)^\gamma}{(\lambda^{n_T}+\vert x\vert)^\gamma}  w_\gamma(x)\, dx\geq \lambda^{n_T(\gamma-1)} \|\vu_n(\lambda^{2n_T}t,.)\|_{L^2_{‘w_\gamma}}^2.$$
 and
 \begin{equation*} \begin{split}   \int_0^{\lambda^{-2n_T}T} \|\vN\vw_n\|_{L^2_{w_\gamma}}^2\, ds =& \int_0^T\int \vert\vN\vu_n(s,x)\vert^2 \lambda^{n_T(\gamma-1)} \frac{(1+\vert x\vert)^\gamma}{(\lambda^{n_T}+\vert x\vert)^\gamma}  w_\gamma(x)\, dx \, ds \\ \geq&  \lambda^{n_T(\gamma-1)}  \int_0^T\| \vN\vu_n \|_{L^2_{‘w_\gamma}}^2\, ds.
 \end{split}\end{equation*} 
 
  Thus, we have a uniform control of $\vu_n$ and of $\vN\vu_n$ on $(0,T)$ for $n\geq n_T$. We may then apply the Rellich lemma (Lemma \ref{spacetime}) and Theorem \ref{stability} to find a subsequence $\vu_{n_k} $ that converges to a global solution of the Navier--Stokes equations. Theorem \ref{weightedNS} is proven.\Endproof
\section{Solutions  of the    advection-diffusion problem with initial data in  $L^2_{w_\gamma}$.}
The proof of Theorem \ref{weightedNS} on the Navier--Stokes problem can be easily adapted to the case of the advection-diffusion problem :

  \begin{theorem}\label{advection}
  Let $0<\gamma\leq 2$. Let $0<T<+\infty$. Let  $\vu_{0}$ be a divergence-free vector field such that $\vu_0\in L^2_{w_\gamma}(\mathbb{R}^3)$ and   $\mathbb{F}$ be  a tensor $\mathbb{F}(t,x)=\left(F_{i,j}(t,x)\right)_{1\leq i,j\leq 3}$ such that $\mathbb{F}\in L^2((0,T), L^2_{w_\gamma})$.  Let $\vb$ be a time-dependent divergence free vector-field ($\vN\cdot\vb=0$) such that $\vb\in L^3((0,T),L^3_{w_{3\gamma/2}})$.
  
Then the   advection-diffusion problem 
 \begin{equation*}  (AD) \left\{ \begin{matrix} \partial_t \vu= \Delta \vu  -(\vb\cdot \vN)\vu- \vN p +\vN\cdot \mathbb{F} \cr \cr \vN\cdot \vu=0,    \phantom{space space} \vu(0,.)=\vu_0
 \end{matrix}\right.\end{equation*}
has a solution $\vu$ such that :
 \begin{itemize} 
 \item[$\bullet$]  $\vu$ belongs to $L^\infty((0,T), L^2_{w_\gamma})$ and $\vN\vu$ belongs to $L^2((0,T),L^2_{w_\gamma})$
 \item[$\bullet$] the pressure $p$  is related to $\vu$, $\vb$  and $\mathbb{F}$ through the Riesz transforms $R_i =\frac{\partial_i}{\sqrt{-\Delta}}$ by the formula
 $$ p=\sum_{i=1}^3\sum_{j=1}^3 R_iR_j(b_iu_j-F_{i,j})$$ 
  \item[$\bullet$] the map $t\in [0,T)\mapsto \vu(t,.)$ is weakly continuous from $[0,T)$ to $L^2_{w_\gamma}$, and is strongly continuous at $t=0$ :
 $$ \lim_{t\rightarrow 0} \|\vu(t,.)-\vu_0\|_{L^2_{w_\gamma}}=0.$$
 \item[$\bullet$]  there exists a non-negative locally finite measure  $\mu$ on $(0,T)\times\mathbb{R}^3$ such that
 \begin{equation*} \partial_t(\frac {\vert\vu\vert^2}2)=\Delta(\frac {\vert\vu\vert^2}2)-\vert\vN \vu\vert^2- \vN\cdot\left( \frac{\vert\vu\vert^2}2\vb\right)-\vN\cdot(p\vu) + \vu\cdot(\vN\cdot\mathbb{F})-\mu.\end{equation*}
 \end{itemize}
  \end{theorem}
 
 \Proof
 Again,  we define $\phi_R(x)=\phi(\frac x R)$, $\vu_{0,R}=\mathbb{P}(\phi_R \vu_0)$ and $\mathbb{F}_R=\phi_R \mathbb{F}$.  Moreover, we define $\vb_R=\mathbb{P}(\phi_R \vb)$. We then solve the mollified problem
 \begin{equation*}  (AD_{R,\epsilon}) \left\{ \begin{matrix} \partial_t \vu_{R,\epsilon}= \Delta \vu_{R,\epsilon}  -((\vb_R*\theta_\epsilon)\cdot \vN)\vu_{R,\epsilon}- \vN p_{R,\epsilon} +\vN\cdot \mathbb{F}_{R,\epsilon} \cr \cr \vN\cdot \vu_{R,\epsilon}=0,    \phantom{space space} \vu_{R,\epsilon}(0,.)=\vu_{0,R}
 \end{matrix}\right.\end{equation*} for which we easily find  a unique solution $\vu_{R,\epsilon}$ in $L^\infty((0,T), L^2)\cap L^2((0,T),\dot H^1)$. Moreover, this solution belongs to $\mathcal{C}([0,T), L^2)$.

Again,  $\vu_{R,\epsilon}$ fulfills the assumptions of Theorem  \ref{estimates} :
\begin{itemize}
 \item[$\bullet$]  $\vu_{R,\epsilon}$ belongs to $L^\infty((0,T), L^2_{w_\gamma})$ and $\vN\vu_{R,\epsilon}$ belongs to $L^2((0,T),L^2_{w_\gamma})$  \item[$\bullet$] the map $t\in [0,T)\mapsto \vu_{R, \epsilon}(t,.)$ is weakly continuous from $[0,T)$ to $L^2_{w_\gamma}$, and is strongly continuous at $t=0$ :
 $$ \lim_{t\rightarrow 0} \|\vu_{R, \epsilon}(t,.)-\vu_{0,R}\|_{L^2_{w_\gamma}}=0.$$
 \item[$\bullet$]   on $(0,T)\times\mathbb{R}^3$, $\vu_{R,\epsilon}$ fulfills the energy equality : 
 \begin{equation*} \partial_t(\frac {\vert\vu_{R,\epsilon}\vert^2}2)=\Delta(\frac {\vert\vu_{R,\epsilon}\vert^2}2)-\vert\vN \vu_{R,\epsilon}\vert^2- \vN\cdot\left( \frac{\vert\vu\vert^2}2\vb_{R,\epsilon}\right)-\vN\cdot(p_{R,\epsilon}\vu_{R,\epsilon}) + \vu_{R,\epsilon}\cdot(\vN\cdot\mathbb{F}_R).\end{equation*}
with $\vb_{R,\epsilon}=\vb_R*\theta_\epsilon$. \end{itemize}
Thus, by Corollary \ref{passive} we know that,
 $$ \sup_{0<t<T} \|\vu_{R,\epsilon}  \|_{L^2_{w_\gamma}} \leq (\|\vu_{0,R}\|_{L^2_{w_\gamma}}+ C_\gamma \|\mathbb{F}_R\|_{L^2((0,T), L^2_{w_\gamma})}) \ e^{C_\gamma (T+ T^{1/3} \|\vb_{R,\epsilon} \|_{L^3((0,T), L^3_{w_{3\gamma / 2}})}^2)}$$
and
 $$  \|\vN\vu_{R,\epsilon}  \|_{L^2((0,T),L^2_{w_\gamma} ) } \leq (\|\vu_{0,R}\|_{L^2_{w_\gamma} }+ C_\gamma \|\mathbb{F}_R \| _{L^2((0,T), L^2_{w_\gamma})}) \ e^{C_\gamma (T+ T^{1/3} \|\vb_{R,\epsilon} \|_{L^3((0,T), L^3_{w_{3\gamma /2}})}^2)}$$ where the constant $C_\gamma$ depends only on $\gamma$.
 
Moreover, we have that 
$$  \|\vu_{0,R}\|_{L^2_{w_\gamma}}\leq C_\gamma  \|\vu_0\|_{L^2_{w_\gamma}},  \|\mathbb{F}_R\|_{L^2_{w_\gamma}} \leq \|\mathbb{F}\|_{L^2_{w_\gamma}}  $${ and }  $$
  \|\vb_{R,\epsilon}\|_{L^3((0,T), L^3_{w_{3\gamma/2}} ) } \leq \|\mathcal{M}_{\vb_R}\|_{L^3((0,T), L^3_{w_{3\gamma/2}} ) }\leq C_\gamma'       \| \vb\|_{L^3((0,T), L^3_{w_{3\gamma/2}} ) }$$
 
 Let $R_n\rightarrow +\infty$ and $\epsilon_n\rightarrow 0$. Let $\vu_{0,n}=\vu_{0,T_n}$, $\mathbb{F}_n=\mathbb{F}_{R_n}$, $\vb_n=\vb_{R_n,\epsilon_n}$ and $\vu_n=\vu_{R_n,\epsilon_n}$.  We may then apply Theorem \ref{stability}, since  $\vu_{0,n}$ is strongly convergent to $\vu_{0}$ in $L^2_{w_\gamma}$,   $\mathbb{F}_n$ is strongly convergent to $\mathbb{F}$  in $L^2((0,T), L^2_{w_\gamma})$,  and   the sequence $\vb_n$ is  strongly convergent  to $\vb$  in $L^3((0,T),   L^3_{w_{3\gamma/2}})$. Thus  there exists $p$, $\vu$  and an increasing  sequence $(n_k)_{k\in\mathbb{N}}$ with values in $\mathbb{N}$ such that

  \begin{itemize} 
 \item[$\bullet$] $\vu_{n_k}$ converges *-weakly to $\vu$ in $L^\infty((0,T), L^2_{w_\gamma})$, $\vN\vu_{n_k}$ converges weakly to $\vN\vu$ in $L^2((0,T),L^2_{w_\gamma})$
 \item[$\bullet$] $p_{n_k}$ converges weakly to $p$ in $L^{3}((0,T),L^{6/5}_{w_{\frac {6\gamma}5}})+L^{2}((0,T),L^{2}_{w_\gamma})$
  \item[$\bullet$] $\vu_{n_k}$ converges strongly  to $\vu$ in  $L^2_{\rm loc}([0,T)\times\mathbb{R}^3)$.  \end{itemize}
 We then easily finish the proof. \Endproof

\section{Application to the study of $\lambda$-discretely self-similar solutions}
  
  We may now apply our results to the study of $\lambda$-discretely self-similar solutions for the Navier--Stokes equations.
  
\begin{definition} Let  $\vu_0\in L^2_{\rm loc}(\mathbb{R}^3)$.  We say that $\vu_0$  is a $\lambda$-discretely self-similar function
($\lambda$-DSS) if there exists  $\lambda>1$ such that   $\lambda\vu_0(\lambda x)= \vu_0$.

A vector field  $\vu\in L^2_{\rm loc}([0,+\infty)\times\mathbb{R}^3)$    is $\lambda$-DSS if there exists  $\lambda>1$  such that $\lambda \vu(\lambda^2 t,\lambda x)=\vu(t,x)$.

A forcing tensor $\mathbb{F}\in L^2_{\rm loc}([0,+\infty)\times\mathbb{R}^3)$    is $\lambda$-DSS if   there exists $\lambda>1$  such that $\lambda^2 \mathbb{F}(\lambda^2 t,\lambda x)=\mathbb{F}(t,x)$.
 \end{definition} 

We shall speak of self-similarity if $\vu_0$, $\vu$ or $\mathbb{F}$ are $\lambda$-DSS for every $\lambda>1$.

$\  $

\noindent{\bf Examples :}
\begin{itemize}
  \item[$\bullet$]   
  Let $\gamma>1 $ and $\lambda>1$. Then, for two positive constants $A_{\gamma,\lambda}$ and $B_{\gamma,\lambda}$, we have :
if $\vu_0\in L^2_{\rm loc}(\mathbb{R}^3)$ is $\lambda$-DSS, then $\vu_0\in L^2_{w_\gamma}$ and
  $$ A_{\gamma,\lambda} \int_{ 1<\vert x\vert\leq \lambda} \vert \vu_0(x)\vert^2\, dx \leq \int \vert \vu_0(x)\vert^2  w_\gamma(x)\, dx \leq B_{\gamma,\lambda} \int_{ 1<\vert x\vert\leq \lambda} \vert \vu_0(x)\vert^2\, dx  $$
   \item[$\bullet$]  $\vu_0\in L^2_{\rm loc}$ is self-similar if and only if it is of the form $\vu_0=\frac{\vw_0(\frac x{\vert x\vert})}{\vert x\vert}$ with $\vw_0\in L^2(S^2)$.
     \item[$\bullet$] $\mathbb{F}$ belongs to $L^2((0,+\infty),L^2_{w_\gamma})$ with $\gamma>1$ and is self-similar if and only if it is of the form $\mathbb{F}(t,x)=\frac 1 t \mathbb{F}_0(\frac x{\sqrt t})$ with $\int \vert\mathbb{F}_0(x)\vert^2 \frac 1{\vert x\vert}\, dx <+\infty$.
   \end{itemize} 
   
   \Proof 
   \begin{itemize}
     \item[$\bullet$]    If $\vu_0$ is $\lambda$-DSS and if $k\in\mathbb{Z}$ we have
                 $$ \!\!\ \!\!\  \int_{\lambda^k<\vert x\vert<\lambda^{k+1}}   \!\!\ \!\!\  \vert\vu_0(x)\vert^2 w_\gamma(x)\, dx \leq \frac {\lambda^k}{(1+\lambda^k)^\gamma}  \int_{1<\vert x\vert<\lambda}   \!\!\ \!\!\ \vert \vu_0(x)\vert^2 \, dx$$
                 with  $\sum_{k\in\mathbb{Z}}  \frac {\lambda^k}{(1+\lambda^k)^\gamma}<+\infty$ for $\gamma>1$.
          \item[$\bullet$]   If $\vu_0$ is self-similar, we have $\vu_0(x ) =\frac 1{\vert x\vert} \vu_0(\frac x{\vert x\vert})$.    From this equality, we find that, for $\lambda>1$
          $$ \int_{1<\vert x\vert<\lambda} \vert\vu_0(x)\vert^2\, dx =(\lambda-1) \int_{S^2} \vert \vu_0(\sigma) \vert ^2 \, d\sigma $$ 
             \item[$\bullet$]   If  $\mathbb{F}$   is self-similar, then it is of the form $\mathbb{F}(t,x)=\frac 1 t \mathbb{F}_0(\frac x{\sqrt t})$. Moreover, we have
             $$\!\!\ \!\!\   \int_0^{+\infty} \!\!\! \int  \vert \mathbb{F}(t,x)\vert^2\, w_\gamma(x) \, dx\, ds=\int_0^{+\infty} \!\!\ \int \vert \mathbb{F}_0(x)\vert^2 w_\gamma(\sqrt t\,  x) \, dx \, \frac{dt}{\sqrt t} = C_\gamma \int \vert\mathbb{F}_0(x)\vert^2\, \frac {dx}{\vert x\vert}$$ with $ C_\gamma=\int_0^{+\infty} \frac 1{(1+ \sqrt\theta)^\gamma} \frac{d\theta}{\sqrt\theta}<+\infty$.  \Endproof
          \end{itemize}
   
   In this section, we are going to give a new proof of the results of Chae and Wolf \cite{CW18} and Bradshaw and Tsai \cite{BT19} on the existence of $\lambda$-DSS solutions of the Navier--Stokes problem (and of Jia and \v{S}ver\'ak \cite{JS14} for  self-similar solutions) :

 \begin{theorem} \label{selfsimilar} Let $4/3<\gamma\leq 2  $ and $\lambda>1$. If $\vu_{0}$ is a $\lambda$-DSS divergence-free vector field (such that $\vu_0\in L^2_{w_\gamma}(\mathbb{R}^3)$) and if $\mathbb{F}$ is a $\lambda$-DSS tensor $\mathbb{F}(t,x)=\left(F_{i,j}(t,x)\right)_{1\leq i,j\leq 3}$ such that $\mathbb{F}\in L^2_{\rm loc}([0,+\infty)\times\mathbb{R}^3)$ , then the Navier--Stokes equations with initial value $\vu_0$
 \begin{equation*}  (NS) \left\{ \begin{matrix} \partial_t \vu= \Delta \vu  -(\vu\cdot \vN)\vu- \vN p +\vN\cdot \mathbb{F} \cr \cr \vN\cdot \vu=0,    \phantom{space space} \vu(0,.)=\vu_0
 \end{matrix}\right.\end{equation*}
 has a global weak solution $\vu$ such that :
 \begin{itemize} 
 \item $\vu$ is a $\lambda$-DSS vector field
 \item[$\bullet$] for every $0<T<+\infty$, $\vu$ belongs to $L^\infty((0,T), L^2_{w_\gamma})$ and $\vN\vu$ belongs to $L^2((0,T),L^2_{w_\gamma})$
 \item[$\bullet$] the map $t\in [0,+\infty)\mapsto \vu(t,.)$ is weakly continuous from $[0,+\infty)$ to $L^2_{w_\gamma}$, and is strongly continuous at $t=0$ :
 $$ \lim_{t\rightarrow 0} \|\vu(t,.)-\vu_0\|_{L^2_{w_\gamma}}=0.$$
 \item[$\bullet$] the solution $\vu$ is suitable : there exists a non-negative locally finite measure $\mu$ on $(0,+\infty)\times\mathbb{R}^3$ such that
 $$ \partial_t(\frac {\vert\vu\vert^2}2)=\Delta(\frac {\vert\vu\vert^2}2)-\vert\vN \vu\vert^2- \vN\cdot\left( (\frac{\vert\vu\vert^2}2+p)\vu\right) + \vu\cdot(\vN\cdot\mathbb{F})-\mu.$$
 \end{itemize} 
 \end{theorem}
 \subsection{The linear problem.}
 Following Chae and Wolf, we consider  an approximation of the problem that is consistent with the scaling properties of the equations : let $\theta$ be a  non-negative and radially decreasing  function in $\mathcal{D}(\mathbb{R}^3)$ with $\int\theta\, dx=1$; We define $\theta_{\epsilon,t}(x)=\frac 1 {(\epsilon\sqrt t)^3}\ \theta(\frac x{\epsilon\sqrt t})$.  We then will study the ``mollified'' problem
 \begin{equation*}  (NS_\epsilon) \left\{ \begin{matrix} \partial_t \vu_\epsilon= \Delta \vu_\epsilon  -((\vu_\epsilon*\theta_{\epsilon,t}) \cdot \vN)\vu_\epsilon- \vN p_\epsilon +\vN\cdot \mathbb{F} \cr \cr \vN\cdot \vu=0,    \phantom{space space} \vu(0,.)=\vu_0
 \end{matrix}\right.\end{equation*} and begin with the linearized problem
 \begin{equation*}  (LNS_\epsilon) \left\{ \begin{matrix} \partial_t \vv= \Delta \vv  -  ((\vb*\theta_{\epsilon,t})\cdot \vN)\vv- \vN q +\vN\cdot \mathbb{F} \cr \cr \vN\cdot \vv =0,    \phantom{space space} \vv(0,.)=\vu_0
 \end{matrix}\right.\end{equation*}

  \begin{lemma}\label{solvLNS}
  Let $1<\gamma\leq 2$.  Let $\lambda>1$  Let  $\vu_{0}$ be a $\lambda$-DSS divergence-free vector field such that $\vu_0\in L^2_{w_\gamma}(\mathbb{R}^3)$ and   $\mathbb{F}$ be  a $\lambda$-DSS  tensor $\mathbb{F}(t,x)=\left(F_{i,j}(t,x)\right)_{1\leq i,j\leq 3}$ such that, for every $T>0$, $\mathbb{F}\in L^2((0,T), L^2_{w_\gamma})$.  Let $\vb$ be a $\lambda$-DSS  time-dependent divergence free vector-field ($\vN\cdot\vb=0$) such that, for every $T>0$,  $\vb\in L^3((0,T),L^3_{w_{3\gamma/2}})$.
  
Then the   advection-diffusion problem 
 \begin{equation*}  (LNS_\epsilon) \left\{ \begin{matrix} \partial_t \vv= \Delta \vv  -  ((\vb*\theta_{\epsilon,t})\cdot \vN)\vv- \vN q +\vN\cdot \mathbb{F} \cr \cr \vN\cdot \vv =0,    \phantom{space space} \vv(0,.)=\vu_0
 \end{matrix}\right.\end{equation*}
has a unique solution $\vv   $ such that :
 \begin{itemize} 
 \item[$\bullet$]  for every positive $T$,  $\vv   $ belongs to $L^\infty((0,T), L^2_{w_\gamma})$ and $\vN\vv   $ belongs to $L^2((0,T),L^2_{w_\gamma})$
 \item[$\bullet$] the pressure $p$  is related to $\vv   $, $\vb$  and $\mathbb{F}$ through the Riesz transforms $R_i =\frac{\partial_i}{\sqrt{-\Delta}}$ by the formula
 $$ p=\sum_{i=1}^3\sum_{j=1}^3 R_iR_j((b_i*\theta_{\epsilon,t})v_j-F_{i,j})$$ 
  \item[$\bullet$] the map $t\in [0,+\infty)\mapsto \vv   (t,.)$ is weakly continuous from $[0,+\infty)$ to $L^2_{w_\gamma}$, and is strongly continuous at $t=0$ :
 $$ \lim_{t\rightarrow 0} \|\vv   (t,.)-\vu_0\|_{L^2_{w_\gamma}}=0.$$
 \end{itemize}
 This solution $\vv$ is a $\lambda$-DSS vector field.
  \end{lemma}

 \Proof  As we have $\vert \vb(t,.)*\theta_{\epsilon,t}\vert\leq \mathcal{M}_{\vb(t,.)} $ and thus 
 $$ \|\vb(t,.)*\theta_{\epsilon,t}\|_{L^3((0,T), L^3_{w_{3\gamma/2}})}\leq C_\gamma \|\vb\|_{L^3((0,T), L^3_{w_{3\gamma/2}})}$$ we see that we can use Theorem \ref{advection} to get a solution $\vv$ on $(0,T)$.
 
 As clearly  $\vb*\theta_{\epsilon,t}$ belongs   to $L^2_t L^\infty_x(K)$ for every compact subset $K$ of  $(0,T)\times\mathbb{R}^3 $, we can use Corollary \ref{unique} to see that $\vv$ is unique.

 Let $\vw(t,x)=\frac 1 \lambda \vv(\frac t{\lambda^2}, \frac x \lambda)$. As $b*\theta_{\epsilon,t}$ is still $\lambda$-DSS, we see that $\vw$ is solution of $(LNS_\epsilon)$ on  $(0,T)$, so that $\vw=\vv$. This means that $\vv$ is $\lambda$-DSS. \Endproof
 
 \subsection{The mollified Navier--Stokes equations.}
 
 The solution 
 $\vv$  provided by Lemma \ref{solvLNS} belongs to $L^3((0,T), L^3_{w_{3\gamma/2}})$ (as  $\vv   $ belongs to $L^\infty((0,T), L^2_{w_\gamma})$ and $\vN\vv   $ belongs to $L^2((0,T),L^2_{w_\gamma})$). Thus we have a mapping $L_\epsilon : \vb\mapsto \vv$ which is defined from
 $$ X_{T,\gamma}=\{\vb \in L^3((0,T), L^3_{w_{3\gamma/2}})\ /\ \vb \text{ is } \lambda-\text{DSS}\}$$ to $X_{T,\gamma}$ by $L_\epsilon(\vb)=\vv$.
 
 \begin{lemma} For $4/3<\gamma$, $X_{T,\gamma}$ is a Banach space for the equivalent norms $\|\vb\|_{L^3((0,T),L^3_{w_{3\gamma/2}})}$ and $\|\vb\|_{L^3((0,T/{\lambda^2}),\times B(0,\frac 1 \lambda))}$. \end{lemma}
 \Proof  We have
 $$ \int_0^T \int_{B(0,1)} \vert \vb(t,x)\vert^3\, dx\, dt=  \lambda^2 \int_0^{\frac T{\lambda^2}} \int_{B(0,\frac 1 \lambda)}  \vert \vb(t,x)\vert^3\, dx\, dt $$
 and , for $k\in\mathbb{N}$, 
 $$\int_0^T \int_{\lambda^{k-1}<\vert x\vert<\lambda^{k}} \vert \vb(t,x)\vert^3\, dx\, dt=\lambda^{2k}\int_0^{\frac T{\lambda^{2k}}}\int_{\frac 1 \lambda <\vert x\vert<1} \vert \vb(t,x)\vert^3\, dx\, dt.
 $$ We may conclude, since for $\gamma>4/3$ we have $\sum_{k\in\mathbb{N}} \lambda^{k (2-\frac {3\gamma}2)}<+\infty$.
 
 \begin{lemma}\label{compact} For $4/3<\gamma\leq 2$, the mapping $L_\epsilon$ is continuous and compact on $X_{T,\gamma}$.
 \end{lemma} 
 
 \Proof Let $\vb_n$ be a bounded sequence in $X_{T,\gamma}$ and let $\vv_n=L_\epsilon(\vb_n)$. 
 We remark that the sequence $\vb_n(t,.)*\theta_{\epsilon,t}$ is bounded in $X_{T,\gamma}$.
 Thus, by Theorem \ref{estimates} and Corollary \ref{passive},  the sequence $\vv_n$ is bounded in $L^\infty((0,T), L^2_{w_\gamma})$ and $\vN\vv_n$ is bounded in $L^2((0,T),L^2_{w_\gamma})$.
 
 We now use Theorem \ref{stability} and get that  then there exists $q_\infty$, $\vv_\infty$,    ${\bf  B}_\infty$   and an increasing  sequence $(n_k)_{k\in\mathbb{N}}$ with values in $\mathbb{N}$ such that
  \begin{itemize} 
 \item[$\bullet$] $\vv_{n_k}$ converges *-weakly to $\vv_\infty$ in $L^\infty((0,T), L^2_{w_\gamma})$, $\vN\vv_{n_k}$ converges weakly to $\vN\vv_\infty$ in $L^2((0,T),L^2_{w_\gamma})$
 \item[$\bullet$] $\vb_{n_k}*\theta_{\epsilon,t}$ converges weakly to ${\bf B}_\infty$ in $L^3((0,T), L^3_{w_{3\gamma/2}})$, ,  \item[$\bullet$]  the associated pressures $q_{n_k}$ converge weakly to $q_\infty$ in $L^{3}((0,T),L^{6/5}_{w_{\frac {6\gamma}5}})+L^{2}((0,T),L^{2}_{w_\gamma})$
  \item[$\bullet$] $\vv_{n_k}$ converges strongly  to $\vv_\infty$ in  $L^2_{\rm loc}([0,T)\times\mathbb{R}^3)$ : for every  $T_0\in (0,T)$ and every $R>0$, we have
  $$\lim_{k\rightarrow +\infty} \int_0^{T_0} \int_{\vert y\vert<R} \vert \vv_{n_k}(s,y)-\vv_\infty(s,y)\vert^2\, ds\, dy=0.$$
  \end{itemize}
  
  As $\sqrt{w_\gamma}\vv_n$ is bounded in $L^\infty((0,T),L^2)$ and in $L^2((0,T), L^6)$, it is bounded in $L^{10/3}((0,T)\times\mathbb{R}^3)$. The strong convergence of $\vv_{n_k}$ in $L^2_{\rm loc}([0,T)\times\mathbb{R}^3)$ then implies the strong convergence of $\vv_{n_k}$ in $L^3_{\rm loc}((0,T)\times\mathbb{R}^3)$.
  
  Moreover, $\vv_\infty$ is still $\lambda$-DSS (a property that is stable under weak limits).We find that $\vv_\infty\in X_{T,\gamma}$ and that
  $$ \lim_{n_k\rightarrow +\infty} \int_0^{\frac T {\lambda^2}} \int_{B(0,\frac 1 \lambda)} \vert \vv_{n_k}(s,y)-\vv_\infty(s,y)\vert^3\, ds\, dy=0.$$
  This proves that $L_\epsilon$ is compact.
  
  If we assume moreover that $\vb_n$ is convergent to $\vb_\infty$ in $X_{T,\gamma}$, then necessarily we have ${\bf B}_\infty=\vb_\infty*\theta_{\epsilon,t}$, and $\vv_\infty= L_\epsilon(\vb_\infty)$. Thus, the relatively compact sequence $\vv_n$ can have only one limit point; thus it must be convergent. This proves that $L_\epsilon$ is continuous.\Endproof

 \begin{lemma}\label{apriori} Let $4/3<\gamma\leq 2$. If, for some $\mu\in [0,1]$,  $\vv$ is a solution of $\vv=\mu L_\epsilon(\vv)$ then $$\|\vv\|_{X_{T, \gamma}}\leq C_{\vu_0,\mathbb{F},\gamma, T}$$ where the constant $C_{\vu_0,\mathbb{F},\gamma, T}$ depends only on $\vu_0$, $\mathbb{F}$, $\gamma$ and $T$ (but not on $\mu$ nor on $\epsilon$).
 \end{lemma}
 
 \Proof We have $\vv=\mu \vw$; with 
  \begin{equation*}  \left\{ \begin{matrix} \partial_t \vw= \Delta \vw  -((\vv*\theta_{\epsilon,t}) \cdot \vN)\vw- \vN q +\vN\cdot \mathbb{F} \cr \cr \vN\cdot \vw=0,    \phantom{space space} \vw(0,.)=\vu_0
 \end{matrix}\right.\end{equation*}
 
 Multiplying by $\mu$, we find that
  \begin{equation*}  \left\{ \begin{matrix} \partial_t \vv = \Delta \vv   -((\vv*\theta_{\epsilon,t}) \cdot \vN)\vv - \vN (\mu q) +\vN\cdot \mu  \mathbb{F} \cr \cr \vN\cdot \vv =0,    \phantom{space space} \vv (0,.) =\mu \vu_0
 \end{matrix}\right.\end{equation*}
 
 We then use Corollary \ref{active}. We choose $T_0\in (0,T)$ such that
 $$ C_\gamma  \left(1+\|\vu_0\|_{L^2_{w_\gamma}}^2+\int_0^{T_0} \|\mathbb{F}\|_{L^2_{w_\gamma}}^2\, ds\right)^2\, T_0\leq 1.$$ Then, as
 $$ C_\gamma  \left(1+\| \mu \vu_0\|_{L^2_{w_\gamma}}^2+\int_0^{T_0} \ \| \mu \mathbb{F} \|_{L^2_{w_\gamma}}^2\, ds\right)^2\, T_0\leq 1$$
 we know that
 $$ \sup_{0\leq t\leq T_0} \|\ \vv(t,.)\|_{L^2_{w_\gamma}}^2 \leq
 C_\gamma (1 + \mu^2  \|\vu_0\|_{L^2_{w_\gamma}}^2 +\mu^2\int_0^{T_0} \|\mathbb{F}\|_{L^2_{w_\gamma}}^2\, ds
)$$ and 
$$  { \int_0^{T_0} \|\vN\vv\|_{L^2_{w_\gamma}}^2\, ds }\leq  
 C_\gamma (1 +\mu^2 \|\vu_0\|_{L^2_{w_\gamma}}^2 +\mu^2\int_0^{T_0} \|\mathbb{F}\|_{L^2_{w_\gamma}}^2\, ds).$$
 In particular, we have
 $$ \int_0^{T_0} \|\vv\|_{L^3_{w_{3\gamma/2}}}^3\, ds \leq C_\gamma  T_0^{1/4}  (1 + \|\vu_0\|_{L^2_{w_\gamma}}^2 +\int_0^{T_0} \|\mathbb{F}\|_{L^2_{w_\gamma}}^2\, ds)^{\frac 3 2}.
 $$ As $\vv$ is $\lambda$-DSS, we can go back from $T_0$ to $T$.
 \Endproof

 \begin{lemma}\label{schauder} Let $4/3<\gamma\leq 2$. There is at least one solution $\vu_\epsilon$ of the equation $\vu_\epsilon=L_\epsilon(\vu_\epsilon)$.
 \end{lemma}

\Proof Obvious due to the Leray--Schauder principle (and the Schaefer theorem), since $L_\epsilon$ is continuous and compact and since we have  uniform   a priori estimates for the fixed points of $\mu L_\epsilon$ for $0\leq\mu\leq 1$.\Endproof 

\subsection{Proof of Theorem \ref{selfsimilar}.}
We may now finish the proof of Theorem \ref{selfsimilar}. We consider the solutions $\vu_\epsilon$ of $\vu_\epsilon=L_\epsilon(\vu_\epsilon)$.

By Lemma \ref{apriori}, $\vu_\epsilon$ is bounded  in $L^3((0,T), L^3_{w_{3\gamma/2}})$, and so is $\vu_\epsilon*\theta_{\epsilon,t}$. We then know, by Theorem \ref{estimates} and Corollary \ref{passive},  that the familly $\vu_\epsilon$ is bounded in $L^\infty((0,T), L^2_{w_\gamma})$ and $\vN\vu_\epsilon$ is bounded in $L^2((0,T),L^2_{w_\gamma})$.
 
 We now use Theorem \ref{stability} and get that  then there exists $p$, $\vu$,    ${\bf  B}$   and a decreasing  sequence $(\epsilon_k)_{k\in\mathbb{N}}$ (converging to $0$) with values in $(0,+\infty)$ such that
  \begin{itemize} 
 \item[$\bullet$] $\vu_{\epsilon_k}$ converges *-weakly to $\vu$ in $L^\infty((0,T), L^2_{w_\gamma})$, $\vN\vu_{\epsilon_k}$ converges weakly to $\vN\vu$ in $L^2((0,T),L^2_{w_\gamma})$
 \item[$\bullet$] $\vu_{\epsilon_k}*\theta_{\epsilon_k,t}$ converges weakly to ${\bf B}$ in $L^3((0,T), L^3_{w_{3\gamma/2}})$
  \item[$\bullet$]  the associated pressures $p_{\epsilon_k}$ converge weakly to $p$ in $L^{3}((0,T),L^{6/5}_{w_{\frac {6\gamma}5}})+L^{2}((0,T),L^{2}_{w_\gamma})$
  \item[$\bullet$] $\vu_{\epsilon_k}$ converges strongly  to $\vu$ in  $L^2_{\rm loc}([0,T)\times\mathbb{R}^3)$.
  \end{itemize}
  Moreover we easily see that ${\bf B}=\vu$. Indeed, we have  that $\vu*\theta_{\epsilon,t}$ converges strongly in $L^2_{\rm loc}((0,T)\times \mathbb{R}^3)$ as $\epsilon$ goes to $0$ (since it is bounded by $\mathcal{M}_\vu$ and converges, for each fixed $t$, strongly in $L^2_{\rm loc}(\mathbb{R}^3)$); moreover, we have $ \vert (\vu-\vu_\epsilon)*\theta_{\epsilon,t}\vert\leq \mathcal{M}_{\vu-\vu_\epsilon}$, so that the strong convergence of $\vu_{\epsilon_k}$ to $\vu$ is kept by convolution with $\theta_{\epsilon,t}$ as far as we work on compact subsets of $(0,T)\times\mathbb{R}^3 $ (and thus don't allow $t$ to go to $0$).
  
  Thus, Theorem \ref{selfsimilar} is proven.\Endproof


\end{document}